\documentclass{amsart}
\usepackage{amsmath,amssymb,graphicx,amsfonts,amscd,hyperref}
\usepackage{mathrsfs}
\usepackage{mathrsfs}
\usepackage{dsfont}
\usepackage[latin1]{inputenc}
\newcommand{\TT}{{\mathbb{T}}}

 \DeclareMathOperator{\E}{\mathbb{E}}      %expectation
 \DeclareMathOperator{\cov}{Cov}

\DeclareMathOperator{\CJ}{\mathrm{CJ}}

 \newcommand{\dd}{{\mathrm{d}}}            %differential
 \newcommand{\ii}{{\mathrm{i}}}

\renewcommand{\Re}{{\mathfrak{Re }}}
\renewcommand{\Im}{{\mathfrak{Im }}}

\renewcommand{\det}{\mathrm{det}}

\newtheorem{thm}{Theorem}[section]

\newtheorem{lem}[thm]{Lemma}
\newtheorem{prop}[thm]{Proposition}
\newtheorem{rem}[thm]{Remark}
\theoremstyle{definition}

\theoremstyle{remark}

\def \be{\begin{eqnarray*}}
\def \ee{\end{eqnarray*}}
\def \ben{\begin{eqnarray}}
\def \een{\end{eqnarray}}

\def\sn{^{(n)}}

\def\QED{\hfill\vrule height 1.5ex width 1.4ex depth -.1ex \vskip20pt}

\numberwithin{equation}{section}

\begin{document}
\title[Orthogonal  polynomials and circular  ensembles]
{Limit theorems for orthogonal polynomials related to circular ensembles}

\author{J. Najnudel}
\address{Institut de Math\'ematiques de Toulouse,
Universit\'e Paul Sabatier,
118 route de Narbonne,
F-31062 Toulouse Cedex 9,
 France}
 \email{joseph.najnudel@math.univ-toulouse.fr}

\author{A. Nikeghbali}
\address{Institut f\"ur Mathematik,
 Universit\"at Z\"urich, Winterthurerstrasse 190,
 CH-8057 Z\"urich,
 Switzerland}
 \email{ashkan.nikeghbali@math.uzh.ch}

\author{A. Rouault}
 \address{Universit\'e  Versailles-Saint Quentin, LMV,
 B\^atiment Fermat, 45 avenue des Etats-Unis,
78035 Versailles Cedex}
 \email{alain.rouault@uvsq.fr}

%\today


\maketitle

\maketitle

\begin{abstract}
For a natural extension of the circular unitary ensemble of order $n$, we study  as $n\rightarrow \infty$ the asymptotic behavior of   the sequence of orthogonal polynomials with respect to the spectral measure. The last term of this sequence is the characteristic polynomial. After taking logarithm and rescaling, we obtain a 
%family indexed by $n$ of
 process indexed by $t \in [0,1]$. We show that it %which   
converges to a deterministic limit, and we describe the fluctuations and the large deviations. 
\end{abstract}

\section{Introduction}

In this paper, we study the asymptotic behavior of a triangular array of complex
 random variables originating  from Random Matrix Theory, and more specifically from some models 
of unitary ensembles.  Our results can be viewed  as  an extension of some earlier works on 
the characteristic polynomial of unitary matrices. Indeed, for fixed $n$, our random variables consist of the sequence of monic orthogonal polynomials with respect 
to a random measure associated with a random unitary matrix of size $n$, when 
evaluated at the point $z=1$, the last term being the characteristic polynomial of that matrix.

It was already noticed in \cite{BHNY} that the characteristic polynomial of a random 
unitary matrix sampled from the Haar measure on $\mathbb U (n)$, the unitary group 
of order $n$ has the same law as a product of independent random variables. In a
 previous paper (\cite{apa}) we saw that this  characteristic polynomial
 $\Phi_n(z)$ (evaluated at $z=1$) is actually a product of variables 
\[\Phi_n(1) = \prod_{j=0}^{n-1} (1-\gamma_j)\]
and  we named the independent variables $\gamma_j$ the modified Verblunsky coefficients, referring to the coefficients involved in the Schur recursion for orthogonal polynomials on the unit circle (OPUC), as dubbed by Simon \cite{Simon1}. 
These variables (for $j < n$) have an explicit  density in the open unit disc, depending on their rank $j$.   We proved in \cite{apa} that this description stays valid when we change the probability measure  on $\mathbb U (n)$, by considering the Circular Jacobi Ensemble, that we will define below.  The construction of the deformed Verblunsky coefficients uses the whole sequence $\{\Phi_{k,n} (z), k=0, \dots, n\}$ of monic orthogonal polynomials with respect to the spectral measure of the pair $(U, e_1)$ where $e_1$ is a fixed vector. 
We have (\cite{apa} Prop. 2.2)
\[\Phi_{k,n}(1) =  \prod_{j=0}^{k-1} (1-\gamma_j)\]
It is then natural to study the triangular array  $\{\Phi_{k,n} (1), k=0, \dots, n\}$ and  in order to normalize the time for different values of $n$,   
 one can study the {\it process} $\{\Phi_{\lfloor nt\rfloor, n} (1) , t \in [0,1]\}$, where $\lfloor \; \rfloor$ denotes the integer part, and 
where, by convention, $\Phi_{0,n} (1) = 1$.

 One can note that $\Phi_{n,n}$ is the characteristic polynomial of $U$ and that when $U$ is the chosen according to the Haar measure, the sequence of random variables $(\Phi_{n,n}(1))$ has played a crucial role in the recent interactions between random matrix theory and number theory (see \cite{KeatingS}). It is also worthwhile to  note that in \cite{killipclock}, Killip and Stoiciu considered also a stochastic process indexed by $t = k/n$ and related to the sequence of orthogonal polynomials in the C$\beta$E ensemble. In fact they considered variables which
are the complex conjugate of our deformed Verblunsky coefficients as auxiliary
variables in the study of the Prüfer phase (Lemma 2.1 in \cite{killipclock}). This tool is used again in \cite{killipfluc}. Later Ryckman (\cite{Ryckman1} Section 4) used a version of the deformed Verblunsky coefficients in the proof of its joint asymptotic laws.

To be more precise and to explain the interest of our approach, we  now  describe our model. 
For $n \geq 1$, $\beta > 0$ and $\delta\in \mathbb C$ such that $\Re\!\ \delta > -1/2$, we consider a
 distribution $\CJ_{\beta, \delta}\sn$ 
on the set of probability measures on the unit circle $\mathbb{T}$, supported
 by $n$ points. This family of distribution generalizes the
 Circular Jacobi Ensemble (the notation CJ comes from this fact), and it can be defined as follows. If the 
random measure 
$$\mu = \sum_{j=1}^n \pi_j \, \delta_{e^{i \theta_j}},$$
for $\theta_j \in [0, 2\pi)$, has the distribution CJ$_{\beta, \delta}\sn$, then:
\begin{itemize}
\item The joint density $h_{\delta,\beta }^{(n)}$ of the law of $(\theta_1, \dots, \theta_n)$, with
 respect to the Lebesgue measure on 
$[0,2 \pi)^n$, is given by 
\begin{equation}\label{loiHP}
h_{\delta,\beta }^{(n)}(\theta_1, \dots, \theta_n) = c_{\delta,\beta
}^{(n)}|\Delta(e^{\ii\theta_1},\dots,
e^{\ii\theta_n})|^\beta\prod_{j=1}^{n}(1-e^{-\ii\theta_j})^
{\delta}(1-e^{\ii\theta_j})^{\overline{\delta}},
\end{equation}
where $c_{\delta,\beta}^{(n)} > 0$ is a  normalization constant. 
\item The weights $(\pi_1, \dots, \pi_n)$ follow a Dirichlet law of parameter $\beta'$ on the simplex $\pi_1 + \dots + \pi_n = 1$, $\pi_j > 0$, 
where $\beta' := \beta/2$ (we conserve this notation $\beta'$ in all the sequel of the paper).  
\item The tuples $(\theta_1, \dots, \theta_n)$ and $(\pi_1, \dots, \pi_n)$ are independent. 
\end{itemize}
Besides, when $U$ is an unitary matrix and $e_1$ a cyclic vector for $U$, then
 the spectral measure $\mu$  of the pair $(U, e_1)$ is defined as the unique (probability) measure on $\mathbb T$ such that 
\[\langle e_1 , U^j e_1\rangle = \int_\mathbb T z^j \mu(dz) \ \ j \in \mathbb Z\,.\]
When $\delta = 0$, it was proved in \cite{killip1} that the distribution of the spectral measure of the pair $(U, e_1)$ where $U$ is randomly sampled from $\mathbb U (n)$ according to the Haar measure and $e_1$ is a fixed vector of $\mathbb C^n$, for instance $(1,0, \dots, 0)$, is precisely  $\CJ_{2, 0}\sn$. 
For $\delta =0$ and $\beta > 0$, these authors found a model of random unitary matrices such that the spectral measure has the distribution $\CJ_{\beta, 0}\sn$. In \cite{apa}, we gave a model corresponding to the case $\delta \not=0$.  All these constructions rely on the theory of orthogonal polynomials on the unit circle (OPUC) that we recall now. 

From the linearly independent family of monomials
$\{1, z, z^2, \dots, z^{n-1}\}$ in $L^2(\mathbb T, \mu)$,
we construct an orthogonal basis $\Phi_{0,n}, \dots, \Phi_{n-1, n}$ of
monic polynomials by the Gram-Schmidt procedure.
The $n^{th}$ degree polynomial obtained this way is
$$\Phi_n (z) =\Phi_{n,n}(z) = \prod_{j=1}^n (z - e^{i\theta_j}),$$
i.e. the characteristic polynomial of $U$. The $\Phi_k$'s
($k= 0, \dots, n$)
obey the Szeg\"o (or Schur) recursion relation:
\be
%\label{Szego}
\Phi_{j+1,n}(z) = z\Phi_{j,n}(z) - \bar\alpha_j \Phi_{j,n}^*(z)
\ee
where \be
%\label{etoile}
\Phi_{j,n}^*(z) = z^j\!\ \overline{\Phi_{j,n}(\bar z^{-1})}\,.
\ee
The coefficients $\alpha_j$ ($j \geq 0)$ are called   Schur
or Verblunsky coefficients and satisfy the condition
$\alpha_0, \cdots , \alpha_{n-2} \in \mathbb D := \{ z\in \mathbb C : |z| < 1\}$ and
$\alpha_{n-1} \in \mathbb T$.
There is a bijection between this set of coefficients and the
set of spectral probability measures $\nu$ (Verblunsky's theorem).
We can write the orthogonal polynomials with the help of a new system of functions built from the Verblunsky coefficients.
Setting 
\[y_k(z) = z - \frac{\Phi_{k+1,n}(z)}{\Phi_{k,n}(z)} = \bar{\alpha}_k\frac{\Phi_{k,n}^*(z)}{\Phi_{k,n}(z)}\ , \ \ (k= 0, \dots, n-1)\,,\]
we have $y_0(z) = \bar\alpha_0$ and the following decomposition:
\be
\Phi_{k,n} (z) = \prod_{j=0}^{k-1} (z - y_j(z))\ , \ \ k=1, \dots, n\,. 
\ee
If $\gamma_j := y_j(1)$, %and $k = n-1$,
 we get
%In particular, for $k = n-1$ we get
\ben
\label{masterf}
\Phi_{k,n}(1) =  \prod_{j=0}^{k-1} (1 - \gamma_j)\ , \ k=1, \dots n\een
and in particular
\ben
\label{deti-u}
\det(I-U) = \Phi_n (1) = \Phi_{n,n}(1)=  \prod_{j=0}^{n-1} (1 - \gamma_j)\,.
\een
%\det(I-U) &=& \prod_{j=0}^{n-1} (1 - y_j)
%\een
%where
%The $y$'s may be also defined by means of the $\alpha$'s with the recursion :
%\ben
%y_k(z) &=& \bar\alpha_k\prod_{j=0}^{k-1}\frac{1-zy^*_j(z)}{z-y_j(z)}\\
%\label{wt} y^*_k(z) &=& \overline{y_k({\bar z}^{-1})}
%\,. \een 
Note that the definition  %(\ref{etoile})
implies
$|\gamma_k|=|\alpha_k|$, and in particular
%In particular, if $y_j := y_j(1)$ and $k = n-1$, we get
%In particular, for $k = n-1$ we get
%\ben
%\label{masterf}
%\det(I-U) = \prod_{j=0}^{n-1} (1 - y_j)
%\een
%where
%\ben
%\label{yalpha}
 %y_0 =\bar\alpha_0\ \ , \ \ y_j = \bar\alpha_j e^{\ii \varphi_{j-1}}\ \ , \ \ e^{\ii \varphi_{j-1}} = \prod_{r=0}^{j-1}\frac{1-{\bar y}_r}{1-y_r} \ \ , \ \ (j=1, \dots, n-1)\,.
%\een
%The last term is particular. Since
% $|\alpha_{n-1}| = 1$,  we have also
 $|\gamma_{n-1}|= 1$.
%put 
%$\alpha_{n-1}= e^{\ii \psi_{n-1}}$, 
%so that \ben\label{last} y_{n-1} = e^{\ii (-\psi_{n-1}+ \varphi_{n-2}%(\alpha_0, \dots, \alpha_{N-2})
%)} := e^{\ii \theta_{n-1}}\,.\een
In the sequel, following \cite{apa}, we refer to the $\gamma_j$'s as the \textit {deformed Verblunsky coefficients}.

 In \cite{apa}, it is proven that for $\mu$ following the distribution CJ$_{\beta, \delta}\sn$, the 
coefficients $(\gamma_j)_{0 \leq j \leq n-1}$ 
are independent (note that in general it is not true for the classical Verblunsky coefficients, except if $\delta = 0$) and their distributions %of $\gamma_j$ is 
are explicitly computable. More precisely,  for $r > 0$ and $\delta \in \mathbb C$ such that $r+2\Re\!\ \delta + 1 > 0$, let $g_r^{(\delta)}$ be the density on the unit disc $\mathbb D$ proportional to
\[\left(1-|z|^2\right)^{r-1} (1-z)^{\bar\delta} (1-\bar z)^{\delta}\] 
and let $\lambda^{(\delta)}$ be the density on the unit circle $\mathbb U$ proportional to
\[ (1-z)^{\bar\delta} (1-\bar z)^{\delta}\,.\]
Then, for $j < n-1$, $\gamma_j$ has density $g_{\beta' (n-j-1)}^{(\delta)}$ and $\gamma_{n-1}$ has density $\lambda^{(\delta)}$. Note that the law of $\gamma_{n-1}$ does not depend on $n$ and $\beta$. The values of the normalization constants
 can 
 be easily deduced from the computation of integrals which are collected in the appendix of the present paper. 

Let us now explain why we choose to focus on the value at $1$ of the characteristic polynomial (in addition to its major role in the number theory connections mentioned earlier). Note
 that in the case $\delta = 0$ (which corresponds to the classical 
Circular Jacobi Ensemble, and in particular the Circular Unitary Ensemble for $\beta =2$), the law of $\Phi_n(z)$ does not
depend on $z \in \mathbb{T}$, since the distribution of the eigenvalues of $U$ is invariant by rotation. 
On the other hand, for $\delta \neq 0$, the point $1$ plays an important 
role since it is a singularity of the potential.  It is then classical to study the behavior of $\Phi_{n}(1)$ in the large $n$ asymptotics (see for example \cite{KeatingS}, \cite{HKOC}, \cite{killipfluc}, \cite{Ryckman1}). Notice that all these authors consider the case $\delta = 0$. 

Our extension to a study of the array $\{\Phi_{k,n}(1), k \leq n\}$ has its own interest as a study of characteristic polynomials. It comes from the following remark.  
From a measure $\mu$ carried by $n$ points, one can also define a $n \times n$ unitary matrix $U_n$, called GGT by Simon \cite{simon2006cmf} section 10, and which is the matrix of 
the linear application $h$ on $L^2(\mathbb{T}, \mu)$ given by $h(f)(z) = z f(z)$, taken in a basis of orthonormal polynomials with respect to $\mu$. 
If $1 \leq k \leq n$, one can denote by $G_k(U_n)$ the $k \times k$ topleft submatrix of $U$ (which is not unitary in general).
%, and by $\Phi_k$ its characteristic polynomial:
%$$\Phi_{k,n}(z) = \operatorname{det} (I_k - z G_k(U)).$$
Then  it is known (see for example Forrester \cite{Forrester}  Prop. 2.8.2, or Simon \cite{simon2006cmf} proof of Prop. 3.1)
that one has for all $k \in \{1,2,\dots, n\}$:
\ben
\label{detprody}
\Phi_{k,n} (1) = \det(I_k - G_k(U_n)) %= \prod_{j=0}^{k-1}(1-\gamma_j),
\een
%and in particular, 
%\Phi_{n,n} (1) = \det(I_n - U) = \prod_{j=0}^{n-1}(1-\gamma_j) =: \Phi_n (1)\,,$$
%where the $\gamma_j$'s are the \textsl{modified} Verblunsky coefficients.  
For other aspects of this model see \cite{BNR}.

Most of the results we will obtain on the process  $\{\Phi_{\lfloor nt\rfloor, n} (1) , t \in [0,1]\}$ will in fact concern its logarithm. Note 
that even when $\Phi_{k,n}(1) \neq 0$, its complex logarithm is not obvious to define rigorously, since its imaginary part is a priori given only modulo $2 \pi$.
However, there is a natural way to deal with this issue, which is described in the appendix. One can then fully justify the following formula:
\begin{equation}
\label{logdetprody}
\log \Phi_{k,n} (1) = \sum_{j=0}^{k-1} \log(1-\gamma_j),
\end{equation}
when $\Phi_{k,n}(1) \neq 1$, which occurs almost surely under CJ$_{\beta, \delta}\sn$. 

 We study asymptotic properties of this determinant as $n \rightarrow \infty$ under essentially two regimes:
 
$\bullet$ First regime: $\delta$ is fixed and $\Re\!\  \delta > -1/2$ (hence, this regime includes the case 
$\delta = 0$). 

$\bullet$ Second regime: $\delta = \beta' \dd n$ with $\Re\!\ \dd > 0$.

Some of the  results proved here were already announced in \cite{notcras}.  When $\beta= 1,2,4$, the independence of the random variables $\Phi_{k+1, n}/ \Phi_{k,n}$ $k=0, \dots , n-1$ and the identification of their distributions are strongly related to the results of Neretin (\cite{Ner1} Corollary 2.1). In that framework, they can be carried easily on models of matrix balls via his Proposition 2.3. Actually our results may be extended to all models where a remarkable separation of variables occur (see \cite{Ner2}). More precisely the paper is organized as follows.

In Section 2 we study the variables $\log \Phi_{\lfloor nt\rfloor, n}(1)$ for both regimes and prove that their expectations converge to some explicit deterministic function of $t$, and then we study the fluctuations of $(\log \Phi_{\lfloor nt\rfloor, n}(1))$ as a stochastic process on the space of c\`adl\`ag $\mathbb R^2$-valued functions.   It appears that in the first regime, one has to distinguish between the case $0<t<1$ and $t=1$; at $t=1$ some  transition is occuring since the normalization is changed. Indeed, anticipating the notation of next section where we write $\log \Phi_{\lfloor nt\rfloor, n} (1) - \mathbb E \log \Phi_{\lfloor nt\rfloor,n}
 (1)= \xi_n (t) + \ii \eta_n(t)$ and $\zeta_n(t) = \begin{pmatrix}\xi_n(t)\\ \eta_n(t)\end{pmatrix}$, we shall prove that in the first regime, $\{\zeta_n(t);\;t<1\}$ converges to some explicit %complex
 Gaussian diffusion $\{\zeta_\dd^0(t);\;t<1\}$ while $\zeta_n(1)/\sqrt{n}$ converges to a %complex
 Gaussian random variable which is independent of the diffusion $\{\zeta_\dd^0(t);\;t<1\}$. In the second regime, $\{\zeta_n(t);\;0\leq t\leq1\}$ converges to some explicit Gaussian diffusion and there is no normalization to perform.

Section 3 is devoted to the establishing Large Deviation Principle (LDP) for the distributions of the real and imaginary parts of $\log \Phi_{\lfloor nt\rfloor, n}(1)$ as a two dimensional random vector with values in the Skorokhod space endowed with the weak topology. We focus our study in the first regime on the case $\delta=0$. Again in the first regime  the case $t=1$ is playing a special role and is not included. From the contraction principle we are then able to deduce LDP for the marginals at fixed time. Our approach is standard: we first compute the normalized cumulant generating function, compute its limit as well as its dual transform, and end with  exponential tightness.

In Section 4 we discuss the connections between the results of Section 3 and existing results on the LDP for the empirical spectral distribution for the circular Jacobi ensemble. 

Section 5 gathers in appendix some properties of remarkable functions  and densities used  in the proofs, in particular the Gamma function $\Gamma$ and the Digamma function $\Psi$.   For $z \in \mathbb{C} \backslash \mathbb{R}_-$, we take the notation: 
% the Digamma function $\Psi = \Gamma'/\Gamma$ and the notation 
\[\ell(z) := \log \Gamma (z)\ ; \ \Psi(z) = \frac{\Gamma'(z)}{\Gamma(z)}\,.\]
The determination of the logarithm in $\ell$ is chosen in the unique way such that $\ell$ is holomorphic 
on $\mathbb{C} \backslash \mathbb{R}_-$ and real on $\mathbb{R}_+^*$. 

All along the paper, we use the entropy function $\mathcal J$ defined by :
%  In the sequel, the following entropy functions will be often used:
\begin{equation}
\label{defJ}
 {\mathcal J}(u) =
     \begin{cases}
         u\log u -u + 1 & \text{if $u > 0$} \\
      1  & \text{if $u=0$}\\
  +\infty   & \text{if $u<0$}
     \end{cases}
\end{equation}
and its primitive \ben \label{defF} F(t) = \int_0^t {\mathcal
J}(u)\ du = \frac{t^2}{2}\log t -\frac{3t^2}{4} + t, \ \ \ (t \geq
0)\,. \een
When the arguments of $\mathcal J$ or of $F$ are complex, we choose the principal determination of the logarithm. 

\section{Convergence and fluctuations}

 Let, for $\Re\!\  \dd > 0$ and $0 < t \leq 1$,  %$h(x) = x \log x$ and
\ben\nonumber
{\mathcal F}_\dd (t) &=& \log(1+2\Re\!\  \dd) -\log(1+\bar\dd) -\log(1-t+ 2\Re\!\ \dd) + \log(1-t+\bar\dd)\\  \label{defed}\\
\nonumber
{\mathcal E}_\dd(t) &=&  {\mathcal J}(1+2 \Re\!\ \dd) - {\mathcal J}(1 + 2\Re\!\ \dd -t) -{\mathcal J}(1+\bar\dd) + {\mathcal J}(1 +\bar \dd-t)\,. 
\een
In this section, we are interested in the process
\[\zeta_n(t) = \begin{pmatrix}\xi_n(t)\\ \eta_n(t)\end{pmatrix}\]
where we have written
\[\log \Phi_{\lfloor nt\rfloor, n} (1) - \mathbb E \log \Phi_{\lfloor nt\rfloor,n}
 (1)= \xi_n (t) + \ii \eta_n(t)\,.\]
%Recall that $\log \Phi_{\lfloor nt\rfloor, n} (1) - \mathbb E \log \Phi_{\lfloor nt\rfloor,n}
% (1)= \xi_n (t) + \ii \eta_n(t)$ and $\zeta_n(t) = \begin{pmatrix}\xi_n(t)\\ \eta_n(t)\end{pmatrix}$. 
As a consequence of the result just below,  $\mathbb E \log \Phi_{\lfloor nt\rfloor,n}
 (1)$ is finite, and then $\zeta_n(t)$ is well-defined.
\subsection{Convergence to a deterministic limit}
\begin{thm}
\label{2.1}
\begin{enumerate}
\item In the first regime, i.e. for fixed $\delta$, and for $n$ going to infinity, 
%\begin{eqnarray}
\ben
\label{2.2}\E \log \Phi_n (1)= \frac{\delta}{\beta'} \log n +  C + o(1)\een
where $C$ is a constant, and for $0 < t < 1$,
\ben
\label{2.3}
\E \log \Phi_{\lfloor nt\rfloor, n}(1) &=& -\frac{\delta}{\beta'}\log (1-t) + o(1)\,. 
\een
\item In the second regime, i.e. $\delta = \beta' \dd n$, and for $0 <  t  \leq 1$, we have
\ben
\lim_{n  \to \infty}\!\ \left(\mathbb E \log  \Phi_{\lfloor nt \rfloor, n }(1)-n{\mathcal E}_\dd
\left(\frac{\lfloor nt\rfloor}{n}\right)\right) = \left(\frac{1}{2} - \frac{1}{\beta}\right) 
{\mathcal F}_\dd (t)\,,
\een
uniformly in $t$.
\item In the first regime, we have 
\ben
\lim_{n \to \infty} \frac{1}{\log n} \!\ \cov\!\ \zeta_n (1) = \frac{1}{\beta}I_2,
\een 
where $I_2$ denotes the $2 \times 2$ identity matrix, and for $0 < t < 1$,
\ben
\lim_{n \to \infty} \!\ \cov\!\ \zeta_n (t) = \int_0^t {\mathcal Z}_s^0 ds,
\een
where $${\mathcal Z}_t^0 := \frac{1}{\beta(1-t)} I_2.$$
\item In the second regime, we have for $0 < t \leq 1$: 
\ben
\lim_{n \to \infty} \!\ \cov\!\ \zeta_n (t) = \int_0^t {\mathcal Z}_s^\dd ds
\een
where
\ben
\label{defZ}{\mathcal Z}_t^\dd =  \frac{1}{\beta'}\left(\begin{array}{cr}\displaystyle \frac{1}{1-t+ 2\Re\!\ \dd}-\Re\!\  \frac{1}{2(1-t+\dd)} & \ \displaystyle\Im\!\ \frac{1}{2(1-t+\dd)}\\ & \\ \displaystyle \Im\!\ \frac{1}{2(1-t+\dd)} & \ \displaystyle\Re\!\  \frac{1}{2(1-t+\dd)}\end{array}\right)\,.\een
\item In the second regime, 
\ben
\label{2.7}
 \sup_{t\in [0,1]}\Big|\frac{1}{n} \log \Phi_{\lfloor nt \rfloor, n}(1) - 
{\mathcal E}_\dd (t)\Big| \underset{n \rightarrow \infty}{\longrightarrow} 0
\een
in $L^2$, and then in probability.
\end{enumerate}
\end{thm}

%\begin{proof}
\proof
Proof of (1) and (2).

 Taking expectations in (\ref{sum}) and  summing up (\ref{m1c}), we have, for $1 \leq m \leq n$, 
%\[\E \log \Phi_{m,n} (1) = \sum_{j=0}^{m-1} \left[
 %\Psi(\beta'(n-j-1)+1 + \delta + \bar\delta) - \Psi((\beta'(n-j-1)+1+\bar\delta) \right].\]
%or
\begin{equation}
\label{psisum}
\mathbb E \log \Phi_{m,n} (1) = \sum_{k = n-m+1}^{n}
 \left[ \Psi\left(\beta' (k-1)+1+  \delta + \bar\delta \right)-  \Psi\left(\beta' (k-1)+1+ \bar \delta \right)\right]\,.
\end{equation}
%\begin{eqnarray}\nonumber\mathbb E \log \Phi_m (1) = \sum_{n-m}^{n-1}
 %\Psi\left(\beta' k+1+ 2a\beta' n\right)-  \Psi\left(\beta' k+1+\bar \dd\beta'  n\right)\\
%\label{psisum}\end{eqnarray}

% The goal is to prove that for $\delta = n\dd\beta'$,
%\ben\label{expd} \mathbb E \log \Phi_{n, \lfloor nt\rfloor}(1) =  n\mathcal{E}_\dd(t) +\left(\frac{1}{2} - \frac{1}{\beta}\right) {\mathcal F}_\dd (t) +o(n)
%\een
%and that for $\delta$ fixed,
%\ben
%\label{exp0}
%\E \log \Phi_n (1) = \frac{\delta}{\beta'}\log n +C + o(1)\\
%\label{expt}
%\E \log \Phi_{n , \lfloor nt\rfloor} (1) = -\frac{\delta}{\beta'}\log (1-t) + o(1)
%\een
%where in (\ref{exp0}), $C$ is a constant and in (\ref{expt}), $t< 1$.
%\bigskip

If in (\ref{psisum}), we keep only the first two terms from the Abel-Plana formula (\ref{AP}), we get,
for $m \leq n-1$, or in the second regime, for $m =n$ large enough,
\begin{eqnarray}\nonumber
\int_{n-m}^n \left[\Psi\left(\beta'(s-1) +1+\delta + \bar\delta\right) - \Psi\left(\beta'(s-1) +1+\bar\delta\right)\right] ds\\
\nonumber
+ \frac{1}{2}\Psi\left(\beta'(n-1)+1+\delta + \bar\delta\right)-  \frac{1}{2}\Psi\left(\beta'(n-m-1)+1+\delta + \bar\delta\right)\\
\label{pp}
-\frac{1}{2}\Psi\left(\beta'(n-1)+1+\bar\delta\right)+  \frac{1}{2}\Psi\left(\beta'(n-m-1)+1+\bar\delta\right)\,.
 \end{eqnarray}
(Note that in the first regime and for $m =n$, we would get the terms $\Psi\left(-\beta'+1+\bar\delta\right)$
 $\Psi\left(-\beta'+1+\delta + \bar\delta\right)$, which are not always well-defined). 
Integrating $\Psi = \ell'$, we get the following expression:
%\begin{eqnarray*} L := \widetilde \ell(\beta'(n-1) + \delta + \bar\delta)&-& \widetilde \ell(\beta'(n-m-1) + \delta+\bar\delta)\\
%&-&\widetilde \ell(\beta'(n-1) + \bar\delta)+ \widetilde \ell(\beta'(n-m-1) + \bar\delta),
%\end{eqnarray*}
\begin{eqnarray} \label{expr}
 L&:=& \widetilde \ell(\beta'(n-1) + \delta + \bar\delta)- \widetilde \ell(\beta'(n-m-1) + \delta+\bar\delta)\\ \nonumber
&&-\widetilde \ell(\beta'(n-1) + \bar\delta)+ \widetilde \ell(\beta'(n-m-1) + \bar\delta),
\end{eqnarray}
%\ben
%\label{expr}
%\een
where $\widetilde \ell(x) := \frac{1}{\beta'}\ell(x+1) + \frac{1}{2}\Psi(x+1)$.
From the Binet formula we have, for $x > -1$,
\ben\label{da}\widetilde \ell(x) = \frac{1}{\beta'} x \log x + \left(\frac{1}{\beta}+\frac{1}{2}\right) \log x - \frac{x-1}{\beta'} + r_1(x)\een
where
\ben
\label{r1}
\  r_1(x) = \frac{1}{\beta'} \int_0^\infty f(s)[e^{-sx} - e^{-s}]\ \! ds + \frac{1}{2}\int_0^\infty e^{-s x}\left(\frac{1}{2}- sf(s)\ \! ds\right)\,.
\een
Now, set
\begin{eqnarray}
\label{i1}I_1(n,m; \delta)&=& {\mathcal I}\left(n-1+\frac{\delta + \bar\delta}{\beta'}\right) - {\mathcal I}\left(n-m-1+\frac{\delta + \bar\delta}{\beta'}\right)\\
\nonumber
&&- {\mathcal I}\left(n-1+\frac{ \bar\delta}{\beta'}\right)+ {\mathcal I}\left(n-m-1+\frac{\bar\delta}{\beta'}\right)
\end{eqnarray}
and
\begin{eqnarray*}
I_2(n,m; \delta) = \left(\frac{1}{\beta}+\frac{1}{2}\right) J_2(n,m; \delta)
\end{eqnarray*}
with
\begin{eqnarray}
\label{j2}
 J_2(n,m; \delta)&=& \log\left (n-1+\frac{\delta + \bar\delta}{\beta'}\right)- \log\left(n-m-1+\frac{\delta + \bar\delta}{\beta'}\right)\\
\nonumber &&-\log\left(n-1+\frac{\bar\delta}{\beta'}\right)+ \log\left(n-m-1+\frac{ \bar\delta}{\beta'}\right).
\end{eqnarray}
We have: 
$$L = I_1(n,m; \delta)  + I_2 (n,m,\delta) + R\,,$$
where 
\begin{eqnarray*}
R &=& r_1(\beta'(n-1) + \delta + \bar\delta)- r_1(\beta'(n-m-1) + \delta+\bar\delta)\\
&&- r_1(\beta'(n-1) + \bar\delta)+ r_1(\beta'(n-m-1) + \bar\delta). 
\end{eqnarray*}
In the sequel, we use several times the trivial estimates (for $c$ fixed and $x$ tending to infinity)
\ben
\label{trivial}
\log (x+c) = \log x + o(1) \  ; \ {\mathcal I} (x+c) ={\mathcal I}(x) + c\log x + o(1)\,.
\een
Let us set $m := \lfloor nt \rfloor$, and $t_n = m/n$, and let us suppose that we are in the first regime. 
If $0 < t < 1$, we have 
\begin{eqnarray*}
I_1(n,m, \delta) & = & \frac{\delta}{\beta'} \log\left( n-1 + \frac{\bar{\delta}}{\beta'} \right)
- \frac{\delta}{\beta'} \log \left( n - n t_n - 1 + \frac{\bar \delta}{\beta'} \right) + o(1) 
\\ \nonumber  & = &- \frac{\delta}{\beta'} \log ( 1 - t) + o(1),
\end{eqnarray*}
and 
$$J_2 (n,m, \delta) = o(1),$$
which implies 
 $$ L = - \frac{\delta}{\beta'} \log ( 1 - t) + R + o(1).$$ 
In the first regime and for $t = 1$, we have to estimate $$\mathbb E \log \Phi_{n} (1) = 
\mathbb E \log \Phi_{n-1,n} (1) +
 \Psi\left(1+  \delta + \bar\delta \right)-  \Psi\left(1+ \bar \delta \right). 
$$
Since the constant $C$ can be modified, it is equivalent to deal with $m = n-1$ or with 
$m=n$. Taking $m = n-1$ gives for some constants $C_1$ and $C_2$, 
\begin{eqnarray*}I_1(n,n-1, \delta) & = & \frac{\delta}{\beta'} \log\left( n-1 + \frac{\bar{\delta}}{\beta'} \right)
+ \mathcal{I} \left(\frac{\delta + \bar \delta}{\beta'} \right) -
\mathcal{I} \left(\frac{\bar \delta}{\beta'} \right)   + o(1) 
\\ \nonumber  & = & \frac{\delta}{\beta'} \log n + C_1 +  o(1),
 \end{eqnarray*}
$$J_2(n,n-1, \delta) = \log\left(\frac{\bar \delta}{\beta'}\right) - 
 \log\left(\frac{\delta + \bar \delta}{\beta'}\right)  + o(1) = C_2 + o(1),$$
which implies, for some constant $C_3$,
$$L = \frac{\delta}{\beta'} \log n + R + C_3 + o(1).$$
Let us now assume that we are in the second regime. For $n$
 large enough, we check the following estimates, uniform in $t \in (0, 1]$: 
$$I_1 (n,nt_n,n \beta' \dd) =  n  {\mathcal E}_\dd(t_n) - J_2(n,n t_n, n \beta' \dd) + o(1),$$
and 
$$J_2(n,nt_n,n \beta' \dd) = {\mathcal F}_\dd (t_n) + o(1),$$
which implies
$$L = n  {\mathcal E}_\dd(t_n) + \left( \frac{1}{\beta} - \frac{1}{2} \right) {\mathcal F}_\dd (t_n) 
+ R + o(1).$$ 
Hence, (1) and (2) are proven, if we check that in any of the previous situations,  
$\mathbb E \log \Phi_{m,n} (1) - L$ and $R$ tend to a constant when $n$ goes to infinity, that this constant 
is zero, except perhaps in the first regime for $t = 1$, and that the convergence is uniform in $t$ in the 
second regime. 
The first quantity can be expressed in function of the last integral term of the Abel-Plana formula,
and the second one comes from the remaining integral term of the Binet formula. 

%{\bf Restent les restes!!}

\underline{First remaining term: Abel-Plana}

We have to prove the convergence of: 
$$ \int_0^\infty \frac{g(n-m+\ii y) - g(n+ \ii y) - g(n-m -\ii y) + g(n-\ii y)}{e^{2\pi y}-1} dy,$$
when $n$ goes to infinity, for
$$g(z) := \Psi(\beta' (z-1) + 1 + \delta + \bar{\delta} ) - 
\Psi(\beta' (z-1) + 1 +\bar{\delta} ),$$
the limit being zero, except in the first regime for $t = 1$. Moreover, in the second regime, 
we want to check that the convergence is uniform with respect to $t$, i.e. with respect 
to $m \in \{1, \dots, n\}$. 
Note that in the second regime, the function $g$ depends on $n$ via $\delta = n \dd \beta'$. 
It is then sufficient to check the following result: for any sequence $(k_n)_{n \geq 1}$ of integers
such that $1 \leq k_n \leq n$,
$$\int_0^\infty \frac{g(k_n+ \ii y) - g(k_n - \ii y)}{ e^{2\pi y}-1} dy $$
is well-defined for all $n \geq 1$, and tends to zero when $n$ goes to infinity, in 
any case and uniformly with respect to the sequence $(k_n)_{n \geq 1}$ 
in the second regime, and when $(k_n)_{n \geq 1}$ tends to infinity
in the first regime. Indeed, we get the desired result by taking successively $k_n= n$ and 
$k_n = n-m$. Note that in the first regime for $t = 1$, we need to take $k_n = 1$ ($m = n-1$), which gives 
an integral independent of $n$, possibly different from zero. 
Now, we have $$g(k_n \pm \ii y) = \Psi(A' \pm \beta' \ii y  )
 - \Psi(A_0 + \ii B_0 \pm \beta' \ii y  ),$$
$A_0, A' > 0$, $B_0 \in \mathbb{R}$ being given by 
$$A_0 + \ii B_0 = \beta' (k_n-1) + 1 +\bar{\delta} $$
and 
$$A' = \beta' (k_n-1) + 1 + \delta + \bar{\delta}.$$
In any case, $A_0$ and $A'$ tend to infinity with $n$, since in the first regime, $k_n$ goes to infinity, 
and in the second regime, it is the case for the real part of $\delta$. Moreover, 
in the second regime, $A$ and $A'$ are greater than $n \beta' \, \Re\!\ \dd$, 
independently from the sequence $(k_n)_{n \geq 1}$. 
Hence, it is sufficient to check that 
$$\int_0^\infty \frac{ \Psi (A + \ii B + \ii C) - \Psi (A + \ii B - \ii C)}{ e^{2 \pi C/ \beta'} - 1}
\, dC$$
is well-defined for $A > 0$ and $B \in \mathbb{R}$, and tends to zero, uniformy in $B$, for 
$A  \rightarrow \infty$. 
From (\ref{Psisum}) we have the following:
%\ben
%\ii\left(\Psi(A+\ii B + 1) - \Psi(A-\ii B+1)\right)= 2B \sum_{k=1}^\infty \frac{1}{(A+k)^2+B^2},
%\een
%and then 
%$$|\left(\Psi(A+\ii B + 1) - \Psi(A-\ii B+1)\right)| \leq 
%so that when $A= \beta'(n-1)+2a\beta' n$ we get
%\ben
%0 < -\ii \int_0^\infty \frac{\Psi(A+\ii \beta' y) - \Psi(A-\ii\beta' y)}{e^{2\pi y}-1} dy \leq \frac{1}{A}\int_0^\infty \frac{2\beta' y}
%{e^{2\pi y}-1}dy
%\een
%which tends to $0$ as $n$ tends to infinity, and idem for the other term with $n-m$.

%Now, for $A,B > 0$ and $C \in \mathbb R$, 
\be
- \ii\left(\Psi(A+\ii B +\ii C ) - \Psi(A+\ii B-\ii C)\right)= 2C \sum_{k=0}^\infty \frac{1}{(A+k+\ii B)^2+C^2}\,.
\ee
To estimate the sum of this series, we first apply the 
crude estimates $|z| \geq |\Im\!\ z|$ and $|z| \geq |\Re\!\ z|$ to the denominator: 
\begin{eqnarray*}
%\label{628}
|(A+k+\ii B)^2+C^2| &\geq&  2 |B| |A+k|\\
%\label{629}
|(A+k+\ii B)^2+C^2| &\geq& |(A+k)^2 - B^2+ C^2|.
\end{eqnarray*}
The first inequality implies, for $A+k \leq  2|B|$, 
$$|(A+k+\ii B)^2+C^2| \geq (A+k)^2,$$
and the second implies, for $A+k > 2 |B|$, 
$$|(A+k+\ii B)^2+C^2| \geq \frac{3}{4} (A+k)^2 + C^2.$$
Hence, 
$$ \left| \sum_{k=0}^\infty \frac{1}{(A+k+\ii B)^2+C^2}  \right|
\leq \frac{4}{3}  \sum_{k=0}^\infty \frac{1}{(A+k)^2} \leq \frac{4}{3} \left(\frac{1}{A^2} + \frac{1}{A} \right).$$
Therefore, 
$$\int_0^\infty \frac{  |\Psi (A + \ii B + \ii C) - \Psi (A + \ii B - \ii C)|}{ e^{2 \pi C/ \beta'} - 1}
\, dC \leq \frac{4}{3}  \left(\frac{1}{A^2} + \frac{1}{A} \right) \int_0^{\infty} \frac{2C} {e^{2 \pi C/ \beta'} - 1} dC,
$$
where the last integral is finite and does not depend on $A$ and $B$. 

%First, 
%assume $|C| \leq A $. Using (\ref{629}) and the classical series-integral comparison, we obtain: 
%\ben
%\sum_{k=1}^\infty  \frac{1}{|(A+k+\ii C)^2+B^2|} \leq \int_0^\infty \frac{dx}{(A+x)^2 - C^2}= 
%\frac{1}{2|C|}\log\frac{A+|C|}{A-|C|}\,. 
%\een
%If $|C| > A $, we split the series in two parts. For %$r \geq 0$ and
% $ k \leq |C| - A + 2$ we use (\ref{628}) and $ |C| (A+k) \geq A^2$ to get
%\be
 %|C| (A+k) \geq A(A-r)
%\ee
%
%It is clear that, for every $k\geq 0$
%\be
%|(A+k+\ii C)^2+B^2| \geq  2 |C| (A+k)\,.
%\ee
%Now, for $r \geq 0$ and $ k \leq |C| - A + r$
%\be
% |C| (A+k) \geq A(A-r)
%\ee
%and
%\be
%\sum_{1 \leq k \leq  |C| - A + 2} \frac{1}{|(A+k+\ii C)^2+B^2|} \leq  \frac{ |C| - A + 2}{2A^2}\,.
%\sum_{k>k_0} \frac{1}{(A+k)^2 - C^2} \leq 2B\int_0
%\ee
%For $k >  |C| - A + 2$, applying (\ref{629}) and again the comparison series-integral, we get
%\ben
%\l%abel{629}
%|(A+k+\ii C)^2+B^2| \geq (A+k)^2 - C^2
%\een
%and %for $k_0 =  \max\{ |C| - A, 1\}$
%\ben
%\sum_{k >  |C| - A + 2}  \frac{1}{|(A+k+\ii C)^2+B^2|}  \leq % 2B \sum_{k>k_0} \frac{1}{(A+k)^2 - C^2}
% 2B\int_{ |C| - A + r -1} \frac{dx}{(A+x)^2 - C^2} =
% \frac{1}{2|C|}\log (2|C|+1)\,.
%\een
%It seems that $r=2$ fits well.

\underline{Second remaining term: Binet}

The terms involved in the expression of $R$ are 
of the form:
\ben
R_1(x,y) = \int_0^\infty f(s) \left[e^{- sx} - e^{-s y}\right]ds
\een
and
\ben
R_2(x,y) = \int_0^\infty \left(\frac{1}{2} - sf(s)\right) \left[e^{- sx} - e^{-s y}\right]ds,\,
\een
where $x= \beta'(n-m-1)  + \alpha$ and $y =  \beta'(n-1)  + \alpha$
with successively $\alpha = \delta + \bar \delta$ and $\alpha = \bar\delta$.
If $t < 1$ or if we are in the second regime, the real parts of $x$ and $y$ tend to infinity with $n$: moreover, 
the convergence is uniform in $t$ in the second regime.
Then, as in \cite{rou},
the dominated convergence theorem allows to conclude that $R_1(x,y)$ and $R_2 (x,y)$ tend to zero, uniformly 
in $t$ in the second regime. 
If $t = 1$ and if we are in the first regime, then $y$ tends to infinity and $x = \alpha$ does not depend on $n$, 
which implies that $R_1(x,y)$ and $R_2(x,y)$ are still converging when $n$ goes to infinity. 
%Let us know consider the case $\delta$ fixed and $m = \lfloor nt\rfloor$ with $t < 1$. It is straightforward to see that the remainig terms tend to zero.
%
%If $m =n$
 %in (\ref{expr}),  we keep only the terms $\widetilde \ell(\beta'(n-1) + \delta + \bar \delta)$ and $\widetilde \ell(\beta'(n-1) + \bar \delta)$
 %(the other terms are constant). Eventually, we get
%$\frac{\delta}{\beta'}\log n$. The remaining terms coming from Abel-Plana or Binet formula have finite limits, as it can be proved using dominated convergence theorem.

Proof of (3) and (4): Computation of covariances.

%The goal is to prove that
%\begin{eqnarray}
%\label{var1}
%\hbox{Var}\!\ \xi_n(t)
%\Re\!\ \log \Phi_{n, \lfloor nt\rfloor}(1)
%&=& -\frac{1}{\beta'}\log \frac{1-t+2\Re\!\ \dd}{1+2\Re\!\ \dd} +  \frac{1}{\beta}\log \frac{|1-t + \dd|}{|1+ \dd|}\\
%\een
%\ben
%\label{var2}
%\hbox{Var}\!\  \eta_n(t)
%\Im\!\ \log \Phi_{n, \lfloor nt\rfloor}(1) 
%&=& -  \frac{1}{\beta}\log \frac{|1-t + \dd|}{|1+ \dd|}\\
%\een
%and
%\begin{eqnarray}
%\label{var3}
%\hbox{Cov}\!\ \left( \xi_n(t), \eta_n(t)
%\Re\!\ \log \Phi_{n, \lfloor nt\rfloor}(1), \Im\!\ \log \Phi_{n, \lfloor nt\rfloor}(1)
%\right) &=& -\frac{1}{\beta}\Im\!\ \left( \log (1-t+ \dd)- \log(1+ \dd)\right)\,. 
%\end{eqnarray}
%We use the independence of the variables $y_j$ and  the Abel-Plana summation again. 
By independence  of the variables $\gamma_j$, we can sum up variances or covariances issued from (\ref{covtheo}).
We can then prove the announced result in the same way as (1) and (2), using 
 the Abel-Plana summation and the Binet formula again, with $\Psi$ replaced by $\Psi'$: we omit the detail. 
%Let us give the details only for the second calculation, following exactly the scheme of the first moment calculus. The other ones are completely similar.
% which yields as in (\ref{expr})
%\[\widetilde\Psi(\left(\beta'(n-1) + \delta\right)\]
%where we have set
%\[\widetilde \Psi(x)=  \frac{1}{\beta'}\Psi(x+1) - \frac{1}{2}\Psi'(x+1)\]
%and from Binet again we have, for $x$ large
%\[\widetilde\Psi(x) = \frac{1}{\beta'}\log x+ \left(\frac{1}{\beta}-\frac{1}{2}\right)\frac{1}{x} + r_2 (x)\,.\]
%\end{proof}

Proof of (5): the process $(\zeta_n(t))_{0 \leq t \leq 1}$ is a two-dimensional martingale. Hence, 
by Doob's inequality:
\begin{eqnarray}
\mathbb{E} \left[ \sup_{t \in [0,1]} | \zeta_n(t) |^2 \right] 
\leq 4 \mathbb{E}[|\zeta_n(1)|^2] &=& 4 \operatorname{Tr} (\cov \zeta_n (1))
\nonumber \\ \underset{n \rightarrow \infty}{\longrightarrow} \frac{1}{\beta'}  \int_{0}^1 \frac{dt}{1 - t + 2 \Re\!\ \dd} 
&= &\frac{1}{\beta'} \log \left( \frac{ 1 + 2 \Re\!\ \dd}{2 \Re\!\ \dd} \right) < \infty. \nonumber
\end{eqnarray}
Moreover, by the uniform convergence (2), 
$$\sup_{t \in [0,1]}  \left|\mathbb E \log  \Phi_{\lfloor nt \rfloor, n }(1)-n{\mathcal E}_\dd
\left(\frac{\lfloor nt\rfloor}{n}\right)\right| 
 \underset{n \rightarrow \infty}{\longrightarrow} \left(\frac{1}{2} - \frac{1}{\beta}\right) 
\sup_{t \in [0,1]} |{\mathcal F}_\dd (t)| < \infty,$$
and one also has
$$\sup_{t \in [0,1]}  n  \left| {\mathcal E}_\dd
\left(\frac{\lfloor nt\rfloor}{n}\right) - {\mathcal E}_\dd (t) \right| 
\leq \sup_{t \in [0,1]} |{\mathcal E}'_\dd(t)| < \infty.$$
Combining these estimates gives the following $L^2$ bound: 
$$\sup_{n \geq 1} \,  \mathbb{E} \left[ \left( \sup_{t \in [0,1]} \left|\log \Phi_{\lfloor nt \rfloor, n} (1) 
- n {\mathcal E}_\dd (t) \right| \right)^2 \right] < \infty,$$
which gives (4) after dividing by $n^2$. 
\subsection{Fluctuations}
We now state a theorem about some limiting distributions related to the process $\{\Phi_{ \lfloor nt\rfloor,n} (1) , t \in [0,1]\}$. Since it can be
 shown by adapting the arguments of Killip and Stoiciu \cite{killipclock} for similar results,
we omit the proof of our theorem, except for the last part.  

Let $D_T$ and $D$ be the space of c\`adl\`ag $\mathbb R^2$-valued functions  on $[0,T]$ and $[0,1)$  respectively, starting from zero.
%in $\mathbb R^2$ vanishing in $0$ and  
The next theorem is about some limiting distributions related to the process $\{\zeta_n(t) , t \in [0,1]\}$ or $\{\zeta_n(t) , t <1\}$.
More precisely let
%\ben
%\label{defZ}{\mathcal Z}_t^\dd =  \frac{1}{\beta'}\left(\begin{array}{cr}\displaystyle \frac{1}{1-t+ 2\Re\!\ \dd}-\Re\!\  \frac{1}{2(1-t+\dd)} & \ \disp%laystyle\Im\!\ \frac{1}{2(1-t+\dd)}\\ & \\ \displaystyle \Im\!\ \frac{1}{2(1-t+\dd)} & \ \displaystyle\Re\!\  \frac{1}{2(1-t+\dd)}\end{array}\right)\,,\een%
 $({\mathcal Z}_t^\dd)^{1/2}$ denote the positive symmetric square root of ${\mathcal Z}_t^\dd$ defined
 in (\ref{defZ}) and let
 ${\mathbf B}_t $ be a standard two dimensional Brownian motion.

\begin{thm}
\label{CLTen1}
\begin{enumerate}
\item
If $\delta = \beta' \dd n$ with $\Re\!\ \dd > 0$ (second regime), then    as $n \rightarrow \infty$
the process $\{ \zeta_n(t)\  ;\ t \in [0,1]\}_n$ converges in distribution  in the Skorokhod space $D_1$ to
 the Gaussian diffusion $\{\zeta_t^\dd\ ;\ t \in [0,1]\}$, solution of the stochastic differential equation: 
\ben
\label{sde1}
d\zeta_t^\dd = ({\mathcal Z}_t^\dd)^{1/2} d{\mathbf B}_t\,.
\een
\item
If $\delta$ is fixed with $\Re\!\ \delta > -1/2$ (first regime), then the joint law of 
the process $\{\zeta_n(t)\  ;\ t < 1\}$ (with trajectories in the Skorokhod space $D$) and the 
variable 
$$ \Theta := \frac{\log \Phi_n (1) - \frac{\delta}{\beta'} \log n}{\sqrt{\log n}} $$
converges, when $n$ goes to infinity, to the joint distribution of $\{\zeta_t^0\ ; \ t< 1\}$ and 
${\mathcal N}(0 ; \beta^{-1})+ \ii\!\ {\mathcal N} (0;\beta^{-1})$, the process and the two gaussian 
variables being independent. 
\end{enumerate}
\end{thm}
Notice that the convergence in law of $\Theta$ is an extension of the celebrated Keating and Snaith result \cite{KeatingS}.

\proof

In order to prove (1), and in (2), the convergence of $\{\zeta_n(t)\  ;\ t < 1\}$ and $\Theta$, taken separately,
we apply a version of the Lindeberg-L\'evy-Lyapunov criterion, available for convergence of 
processes (\cite{JacShi} Chap. 3c). 
For $t< 1$, or in the second regime, it is enough to prove that
\ben
\label{sum41}
 \sum_{k=0}^{\lfloor nt\rfloor-1} \mathbb E |A_k|^4 \underset{n \rightarrow \infty}{\longrightarrow} 0,
\een
where  $$A_k = \log ( 1 - \gamma_k) - \mathbb{E} [\log ( 1 - \gamma_k)].$$
For $t=1$ in the first regime, it is sufficient to check:
\ben
\label{sum42}
\frac{1}{\log^2 n} \sum_{k=0}^{n-1} \mathbb E |A_k|^4 \underset{n \rightarrow \infty}{\longrightarrow} 0.
\een
Now,
\begin{eqnarray*}\mathbb{E} |A_k|^4 &\leq& 8 \mathbb{E} (\Re \, A_k)^4 + 8  \mathbb{E} (\Im \, A_k)^4\\
&\leq& 24  \operatorname{Var}^2 (\Re\!\ \log ( 1 - \gamma_k)) + 24 \operatorname{Var}^2 (
\Re\!\ \log ( 1 - \gamma_k))\\
&&+ 8  \kappa_4  (\Re \log ( 1 - \gamma_k)) +  8 \kappa_4  (\Im \log ( 1 - \gamma_k)), \end{eqnarray*}
where $\kappa_4$ denotes the fourth cumulant. Using the second and fourth order derivatives of 
the function $\Lambda$ introduced in the appendix, we see that $\mathbb{E} |A_k|^4$ is 
a linear combination of terms of the form $\Psi{'''}(r+ 1 + \alpha)$ and $\Psi'(r+ 1 + \alpha) 
\Psi'(r+ 1 + \alpha')  $, 
for $r = \beta'(n-k-1)$ and $\alpha, \alpha' \in \{\delta, \bar \delta, \delta + \bar \delta\}$.  
Now, by (\ref{restepsi}), for $\Re\!\ x > 0$, 
$$\Psi' (x) = \frac{1}{x} + O\left(\frac{1}{(\Re\!\ x)^2}\right)$$
and 
$$\Psi{'''} (x) = \frac{2}{x^3} +  O\left(\frac{1}{(\Re\!\ x)^4}\right)\,.$$
Hence, all the terms involved in the expression of $ \mathbb{E} |A_k|^4 $ are dominated by 
$(n-k)^{-2}$, and they are dominated by $n^{-2}$ in the second regime or for $t < 1$. 
Hence, \[
 \sum_{k=0}^{\lfloor nt\rfloor-1} \mathbb E |A_k|^4
\]
 is dominated by $1/n$, except in the first regime for $t =1$, in which case it is bounded. This shows 
the desired results (\ref{sum41}) and (\ref{sum42}). 
 
In order to prove the convergence of the joint distribution in (2), we can follow the
 scheme of \cite{rou} p.209. We take $0< t_0 < t_1 < 1$. From the above results,
\ben
\frac{\log \Phi_{\lfloor nt_1\rfloor,n }(1)}{\sqrt{\log n}} \rightarrow 0 \label{2806}
\een
in probability, so that
\[\frac{\log \Phi_n (1) - 
\log \Phi_{\lfloor nt_1\rfloor,n}(1)- \frac{\delta}{\beta'}\log n}{\sqrt{\log n}} \rightarrow 
{\mathcal N}(0 ; \beta^{-1}) + \ii {\mathcal N}(0 ; \beta^{-1}),\]
the two gaussian variables being independent.
Now, for $n$ large enough,  
$$\log \Phi_n (1) - \log \Phi_{\lfloor nt_1\rfloor,n}(1) = \sum_{k=\lfloor nt_1\rfloor}^{n-1}( 1-\gamma_k),$$
is independent of $\{\zeta_n(t)\  ;\ t \leq t_0\}$, which is function of $(\gamma_k)_{k < n t_0}$ (recall that
the variables $\gamma_k$ are independent). 
Since, $\{\zeta_n(t)\  ;\ t \leq t_0\}$ tends in law to 
$\{\zeta_t^0\ ; \ t \leq t_0\}$ (as a process with trajectories in $D_{t_0}$), we deduce that 
$$ \left(\{\zeta_n(t)\  ;\ t \leq t_0\} ; \frac{\log \Phi_n (1) - 
\log \Phi_{\lfloor nt_1\rfloor,n}(1)- \frac{\delta}{\beta'}\log n}{\sqrt{\log n}} \right),$$
tends in law to 
$$\left(\{\zeta_t^0\  ;\ t \leq t_0\} ;  {\mathcal N}(0 ; \beta^{-1}) + \ii {\mathcal N}(0 ; \beta^{-1} 
)\right),$$
$\{\zeta_t^0\  ;\ t \leq t_0\}$ being independent of  
${\mathcal N}(0 ; \beta^{-1}) + \ii {\mathcal N}(0 ; \beta^{-1})$.
Using again (\ref{2806}), we deduce that 
$$\left(\{\zeta_n(t)\  ;\ t \leq t_0\} ; \frac{\log \Phi_n (1) -
 \frac{\delta}{\beta'}\log n}{\sqrt{\log n}} \right)$$
has the same limiting distribution. Taking $t_0 \rightarrow 1$ gives  part (2) of the theorem.

%\section{Convergence and fluctuations}

\section{Large deviations}
\subsection{Notation and main statements}
 Throughout this section, we use the standard notation of  \cite{DZ}. In particular we write LDP for Large Deviation
Principle. 
 We say that a sequence $(P_n)_{n\in \mathbb N}$ of probability measures on a Polish space $\mathcal X$ satisfies a LDP with speed  $a_n$ (going to infinity with $n$) and rate function $I$ iff $I: \mathcal X \rightarrow [0, \infty]$ is lower semicontinuous and if 
\[\hbox{For every open set}\ \ O \subset \mathcal X _ , \ \liminf_n \frac{1}{a_n}  \log P_n(O)\geq -\inf_O I\,,\]
\[\hbox{For every closed set}\ \  F \subset \mathcal X _ , \ \liminf_n \frac{1}{a_n}  \log P_n(F)\geq -\inf_F I\,.\]
The rate function is good if its level sets are compact. Moreover, if $X_n$ are $\mathcal X$ random variables distributed according to $P_n$, we say that the sequence $ (X_n)$ satisfies the LDP if the sequence $(P_n)$ satisfies the LDP.

The reader may have some interest in consulting  
\cite{gamb} and mainly \cite{rou} where a similar method is used for a different
model.

For $T \leq 1$, let 
 $\mathrm{M}_T$ stand for the space of signed measures  on $[0,T]$ and let ${\mathrm M}_<$ be the subspace of $\mathrm{M}_1$ consisting  of
measures whose
 support is a compact subset of $[0,1)$.
We endow $D\times D$ with the weak topology $\sigma(D\times D, {\mathrm
M}_<\times {\mathrm
M}_< )$. So, $D\times D$ is the projective limit of the family, indexed by
$T< 1$, of the topological spaces  $\left(D_T\times D_T , \sigma(D_T \times D_T,{\mathrm
M}_T\times  {\mathrm
M}_T)\right)$.

Let $V_\ell$ (resp. $V_r$) be the space of left (resp. right)
continuous $\mathbb{R}$-valued functions with bounded variations. We put a
superscript $T$ to specify the functions on $[0,T]$. There is a
bijective correspondence between $V_r^T$ and ${\mathrm M}_T$:

- for any $v \in V_r^T$, there exists a unique $\mu \in {\mathrm
M}_T$ such that $v = \mu([0, \cdot])$; we denote it by $\dot v$ ,

- for any $\mu \in {\mathrm M}_T$,  $v = \mu([0, \cdot])$ is in
$V_r^T$.

For $v \in D$,  let $\dot v = \dot v_a + \dot v_s$ be the Lebesgue
decomposition of the measure $\dot v$
 in
absolutely continuous and singular parts with respect to the (vectorial) 
Lebesgue measure. The measure  $\dot v_a$ can then be identified with its density.

Now, for $\xi, \eta \in \mathbb{R}$, let us define:
\ben
\label{3.1}
H_a(\xi, \eta) :=  -\xi - \log(2 \cos\eta - e^{\xi})\,,
\een
if $|\eta| < \pi/2$ and $2 \cos\eta - e^{\xi} > 0$, and $H_a(\xi, \eta) = \infty$ otherwise.
% Let also $H_s$ be defined by
%\ben
%H_s(\xi, \eta) = -\xi\,,
%\een
%if $\xi < 0$ and $\eta =0$ and $H_s(\xi, \eta) = \infty$ otherwise.
For $\varphi, \psi \in V_r$ and $T \leq 1$, let us denote:
\begin{equation}
 \label{calh0}
{\mathcal H}_0(T, \varphi, \psi):= \begin{cases}& \int_0^{T} (1-\tau) H_a(\dot\varphi_a (\tau), \dot\psi_a (\tau))\!\ d\tau 
 + \int_0^{T} 
(1-\tau) d (-\dot \varphi_s)(\tau)\\
&\hbox{if} \  \dot \psi_s = 0 \ \hbox{and} \ -\dot \varphi_s  \ \hbox{ is a positive measure}\,,\\
%{\mathcal K}(P,Q) &=  \displaystyle\int_{G}\log\frac{dP}{dQ}dP\;\;\;\;\mbox{if}\;P\ll Q\;\mbox{ and }\log\frac{dP}{dQ}\in L^1(P),\\
      & +\infty\;\;\mbox{ otherwise.}
\end{cases}
\end{equation}
 \begin{thm}
\label{LDPH} 
\begin{enumerate}
\item For every $T < 1$, the sequence %of distributions of
\[\{ n^{-1}\!\ \left(\Re \log \Phi_{\lfloor nt\rfloor,n}(1) , \Im \log \Phi_{\lfloor nt\rfloor,n}(1)\right)\ ; \ t \in [0,T]  \}\] 
under the $C\beta E\sn$ measure (i.e. first regime and $\delta=0$) satisfies 
in $(D_T\times D_T, \sigma(D_T\times D_T, {\mathrm M}_T \times  {\mathrm
M}_T))$
 the LDP  with speed  $\beta' n^2$
with good rate function ${\mathcal H}_0(T, \varphi, \psi)$.
\item For $\Re\!\  \dd > 0$ and $T \leq 1$, the sequence 
% of distributions of
\[\{ n^{-1}\!\ \left(\Re  \log \Phi_{\lfloor nt\rfloor,n}(1) , \Im \log \Phi_{\lfloor nt\rfloor,n}(1)\right)\ t \in [0,T]  \}\] 
under the $\CJ_{\beta, \beta' \dd n}\sn$ measure (second regime) satisfies in
%in $(D, \sigma(D, {\mathrm M}_<))$
 $(D_T\times D_T, \sigma(D_T\times D_T, {\mathrm M}_T \times  {\mathrm
M}_T))$
 the LDP  with speed  $\beta' n^2$
with good rate function%\marginpar{verif delta}
\ben
\label{grf}{\mathcal H}_\dd(T, \varphi, \psi) = {\mathcal H}_0(T, \varphi, \psi) 
-2(\Re\!\  \dd) \varphi(T) -2 (\Im\!\  \dd) \psi(T) -{\mathcal C}_\dd(T)\een
with
\be{\mathcal C}_\dd(T) = 
F(1 + \dd) - F( 1-T+\dd) + F(1 + \bar\dd) - F(1-T+\bar\dd)\\
-F(1 + 2\Re\!\  \dd) + F(1-T+2\Re\!\ \dd) -F(1) + F(1-T)\,.
\ee
\end{enumerate}
\end{thm}

We give now a result on the marginal at time $T$ fixed. It is obtained %To prove the statement of Theorem \ref{mastermarg}, it is enough to
by  applying the contraction principle to the mapping
\[(\varphi, \psi) \mapsto (\varphi(T), \psi(T))\,.\]
In all the sequel of the paper, we consider either the $C\beta E\sn$ ensemble (first regime for $\delta = 0$), 
or the $\CJ_{\beta, \beta' \dd n}\sn$ ensemble for $\Re\!\ \dd > 0$ (second regime). In the first case, we 
put $\dd = 0$. 
\begin{thm}
\label{mastermarg}
 When $(\Re\!\  \dd > 0 , T \leq  1)$ or $(\dd =0, T< 1)$,   the sequence %of distributions of 
\[\{ n^{-1}\!\
 \left(\Re  \log \Phi_{\lfloor nT\rfloor}(1) , \Im  \log \Phi_{\lfloor nT\rfloor}(1)\right)  \}_n\] satisfies the LDP in $\mathbb R^2$ 
 with speed  $\beta' n^{2}$
with good rate function 
\ben\label{proj}h_\dd(T,\xi, \eta) = \inf \{{\mathcal H}_\dd (T, \varphi, \psi) \ | \ \varphi(T) = \xi, \psi(T) = \eta\}\,.\een
In particular (cf. (\ref{grf}))
\ben\label{0d}
h_\dd(T, \xi, \eta) = h_0(T, \xi, \eta) -2(\Re\!\  \dd) \xi -2 (\Im\!\  \dd) \eta -{\mathcal C}_\dd(T), 
\een
which allows to compute easily $h_\dd$ when $h_0$ is known.  
\end{thm}

For the two coordinates, separately, we have the following known result, which
%in which we have corrected a slight mistake in the expression of (1), 
comes from  formula (C.5) in \cite{HKOC}.
\begin{thm}
\label{HKOC}
Assume $\dd = 0$ and $\beta =2$. 
\begin{enumerate}
\item The sequence $\{ n^{-1}\!\
 \Re\!\  \log \Phi_{n}(1)\}_n$
satisfies the LDP in $\mathbb R$ with speed $n^{2}$ and rate function given by the dual (Legendre) of the function:
\[s \mapsto \frac{(1+s)^2}{2} \log(1+s) - \left(1 + \frac{s}{2}\right)^2 \log\left(1+\frac{s}{2}\right) - \frac{s^2}{4} \log \left(2s\right)\]
for $s \geq 0$, and by $\infty$ for $s < 0$. 
It vanishes for negative values of the argument, and it is infinite beyond $\log 2$.
\item The sequence $\{ n^{-1}\!\
 \Im \log \Phi_{n}(1)\}_n$
satisfies the LDP in $\mathbb R$ with speed  $n^{2}$ and good rate function given by the dual (Legendre) of  the function given by
\[t \mapsto \frac{t^2}{8}\log\left(1 + \frac{4}{t^2}\right) - \frac{1}{2}\log\left(1+ \frac{t^2}{4}\right) + t \arctan(t/2)\,.\]
It is finite on $(-\pi/2, \pi/2)$, and infinite otherwise.
\end{enumerate}
\end{thm}

We now give the precise behavior of the first coordinate, in the particular case where $\dd$ is real. 

\begin{thm}
\label{margmoi}
\begin{enumerate}
\item For $(\dd= 0, T < 1)$, or $(\dd > 0, T \leq 1)$, the sequence $\{ n^{-1}\!\
 \Re\!\ \log \Phi_{n}(T)\}_n$
satisfies the LDP in $\mathbb R$ with speed  $\beta 'n^2$ with good rate function $h_\dd(T, \cdot, 0)$.
\item
Let $\xi_T := {\mathcal J}(T) - 1 - {\mathcal J}\left(\frac{1+T}{2}\right) + {\mathcal J}\left(\frac{1-T}{2}\right)\leq 0$.
\begin{enumerate}
\item  If  $\xi \in [\xi_T, T\log 2)$ the equation
\ben\label{implicit}{\mathcal J}(1+\gamma) - {\mathcal J}(1 - T+\gamma) - {\mathcal J}\left(1+\frac{\gamma}{2}\right) + {\mathcal J}\left(1-T+\frac{\gamma}{2}\right) = \xi\een
has a unique solution $\gamma$ and we have
\ben\label{hL1} h_0(T, \xi, 0) = \gamma\xi -{\mathcal L}_0(T, \gamma, 0)\een
where
\ben
\nonumber{\mathcal L}_0(T, \gamma, 0) &:=&F(1+\gamma) - F(1-T+\gamma) - F(1-T) + F(1) \\
\label{hL2}&&-2 F\left(1 + \frac{\gamma}{2}\right) +2F\left(1 -T + \frac{\gamma}{2}\right)\een
\item If $\xi < \xi_T$, then 
\ben\label{lin} h_0(T, \xi, 0) = h_0(T, \xi_T, 0) + (1-T) (\xi_T-\xi)\,.\een
\item If $\xi \geq T\log 2$, then $h_0(T, \xi, 0) = \infty$.
\end{enumerate}
\item For $(\dd > 0, T \leq 1)$ or $(\dd = 0, T < 1)$, the function $\xi \mapsto h_\dd(T, \xi, 0)$ is the dual function of  \[\gamma \mapsto {\mathcal L}_\dd(T, \gamma, 0) := {\mathcal L}_0(T, \gamma +2\dd, 0) - {\mathcal L}_0(T, 2\dd, 0)\,.\]
\end{enumerate}
\end{thm}

\subsection{Proof of Theorem \ref{LDPH}}

We compute the normalized cumulant generating function, find its limit, perform the dual transform, study the exponential tightness and eventually prove that the LDP is satisfied.

%\noindent $\bullet$\textsl{First step: The normalized cumulant generating function}
From (\ref{firstc}) and (\ref{mfy})  we see that 
\ben
\label{RN1}
\frac{\CJ_{\beta, \delta}\left[(1- \bar \gamma_j)^{z} (1- \gamma_j)^{\bar z}\right]}{\CJ_{\beta, 0}\left[(1- \bar \gamma_j)^{(z+\delta)} (1- 
\gamma_j)^{(\bar z + \bar\delta)}\right]}= \frac{c_{r,\delta}}{c_{r, 0}}
\een
 for every $j < n-1$ and $z$ such that %$ n ( 1 + 2\Re\!\ (d + z)) - j -1 > 0$:
 $2\Re\!\ (\delta + z)) > -1$, where 
 $r= \beta'(n-j-1)$. The RHS of (\ref{RN1}) does not depend on $z$. It should then be clear that we can reduce the case $\Re\!\ \dd > 0$ to the case $\dd=0$. The above shift in the argument of the generating function 
provides the linear term $-2(\Re\!\ \dd) \varphi(T) -2 (\Im\!\ \dd) \psi(T)$ in the rate function, and the RHS of (\ref{RN1}) gives the constant $-{\mathcal C}_\dd(T)$.
\smallskip

\noindent $\bullet$ \textsl{First step: The normalized cumulant generating function.}
\begin{lem}Let $T < 1$.
\begin{enumerate}
\item Let us assume $\dd = 0$.
\begin{enumerate}
\item   For every path $(x(\tau), y(\tau))_{\tau \in [0,T]} \in V_\ell^T \times V_\ell^T$ such that $x(\tau) + 1-\tau > 0$ on $[0, T]$, set   $2z(\tau) := x(\tau) + \ii y(\tau)$ and 
\ben
\label{defGL}
&&\ \Lambda_0(T, x, y) :=\\ \nonumber&& \int_0^T \left( {\mathcal J}(1-\tau  + x(\tau)) - {\mathcal J}(1-\tau+ z(\tau)) - {\mathcal J}(1-\tau + \bar z(\tau)) + {\mathcal J}(1-\tau )\right) d\tau\,.
\een
Then we have:
\ben
\label{riemann}\\
\nonumber\lim_{n \rightarrow \infty} \frac{1}{\beta'n^2}\log \mathbb E \exp\left( n \beta' \Re \left[ \int_0^T \left(x(\tau) - \ii y(\tau)\right) d\log \Phi_{\lfloor n\tau\rfloor} (1) \right]
\right)= \Lambda_0(T, x, y)\een 

\item In particular, when $x(.) \equiv s$ and $y(.) \equiv t$,  if we set
%use the notation
 ${\mathcal L}_0(T, s, t) := \Lambda_0(T, x, y)$ (which generalizes the notation ${\mathcal L}_0(T, \gamma, 0)$ introduced in Theorem \ref{margmoi}), 
we have for every $s> -(1-T)$
\ben
\label{lmarg}
{\mathcal L}_0(T, s, t)  &=&F(1 +s) - F(1 -T+s) + F(1) - F(1 -T)
\\ \nonumber
&&- F(1 +z) + F(1 -T+z)- F(1 +  \bar z) + F(1 -T+\bar z)
\een
where $2z = s+\ii t$. %, for every $s> -(1-T)$.
\end{enumerate} 
\item For $\Re\!\  \dd > 0$, the analogues of (\ref{defGL}) and (\ref{lmarg}) are 
\begin{equation}\label{Lambdad}\Lambda_\dd(T,x,y)  = \Lambda_0(T, x+ 2\Re\!\ \dd, y+2\Im\!\ \dd) - \Lambda_0(T, 2\Re\!\ \dd, 2\Im\!\ \dd)\,,\end{equation}
and
\begin{equation}\label{calLambdad}
{\mathcal L}_\dd(T,s,t)  = {\mathcal L}_0(T, s+ 2\Re\!\ \dd, t+2\Im\!\ \dd) - {\mathcal L}_0(T, 2\Re\!\ \dd, 2\Im\!\ \dd)\,,
\end{equation}
for $s > -(1-T) -2\Re\!\ \dd$.
\end{enumerate}
\end{lem}

\proof (1)(a) Let us set $x(\tau) = y(\tau) = z(\tau) := 0$ for $\tau > T$, and  $\tau_j := (j+1)/n$ for $j=0, \dots, n-1$. One has:
\begin{eqnarray*}
&& \mathbb E \exp\left( n \beta' \Re \left[ \int_0^T \left(x(\tau) - \ii y(\tau)\right) d\log \Phi_{\lfloor n\tau\rfloor} (1) \right] \right)
 \\ && = \mathbb{E} \exp \left( n \beta' \Re \left[ \sum_{j=0}^{n-1} \left(x \left(\tau_j\right)  - \ii y (\tau_j) \right) \, \log(1 - \gamma_j) \right] \right)
\\ && = \prod_{j=0}^{n-1} \mathbb{E} \exp \left( n \beta'  \Re \left[ \left(x (\tau_j)  - \ii y (\tau_j) \right) \, \log(1 - \gamma_j) \right] \right)
\end{eqnarray*} 
and using (\ref{core}), we get:
\begin{eqnarray}
&& \log  \mathbb E \exp\left( n \beta' \Re \left[ \int_0^T \left(x(\tau) - \ii y(\tau)\right) d\log \Phi_{\lfloor n\tau\rfloor} (1) \right] \right) \nonumber
\\ && = \sum_{j=0}^{n-1} \left[ \ell( n \beta' (1-\tau_j + x(\tau_j)) + 1) + \ell(n \beta' (1- \tau_j) + 1) \right. 
\nonumber \\ && \left. - \ell( n \beta' (1-\tau_j + z(\tau_j)) + 1 ) -
 \ell(n \beta' (1-\tau_j + \bar{z}(\tau_j)) + 1 ) \right]. \label{nghyi}
\end{eqnarray}
%where $\ell$ denotes the logarithm of the gamma function. 

Now, the Binet formula (\ref{bin1}) yields: %\marginpar{vrai en complexe?}
$$\ell(u+1) =  \log (u) + \log \Gamma(u) = (u+1/2) \log u -u+1 +\int_0^\infty f(s) [e^{-su} - e^{-s}] ds\,.$$ 
%to estimate $\log \mathbb E[R_k^t e^{i\theta_k s}]$.
If we apply four times this formula in order to estimate the term indexed by $j < n$ in \eqref{nghyi}, the contribution of the term $u$ is $0$, the contribution of the term $\log u$ is
\begin{eqnarray*}\log(1-\tau_j)  + \log(1-\tau_j+x(\tau_j))
 -\log(1-\tau_j  +z(\tau_j)) - \log(1-\tau_j+\overline{z(\tau_j)})\,, \end{eqnarray*}
the contribution of the term $u\log u$ is proportional to $\beta'n$ with coefficient
\begin{eqnarray*}
{\mathcal J}(1-\tau_j)  
- {\mathcal J}(1-\tau_j+ z(\tau_j)) - {\mathcal J}(1-\tau_j + \overline{z(\tau_j)}) + {\mathcal J}(1-\tau_j  + x(\tau_j))\,;
\end{eqnarray*}
dividing by $\beta'n^2$ and performing Riemann sums gives the integral in (\ref{defGL}). The remaining
 part is a sum of bounded terms which is negligible with respect to $n^2$. 

(1)(b) The equality (\ref{lmarg}) is obvious by integration.

(2) To get the expression corresponding to $\Re\!\ \dd > 0$ 
we just use (\ref{RN1}).
 \QED
\medskip

\noindent $\bullet$ \textsl{Second step ; the Legendre duality.}

It will be convenient to perform a time-change in (\ref{defGL}), setting
\ben
\label{xtau}
x(\tau)  = (1-\tau)X(\tau) \ ; \ y(\tau)  = (1-\tau)Y(\tau)\ ; \ z(\tau) = (1 - \tau) Z(\tau)\,.
\een
The above expression (\ref{defGL})) of $\Lambda_0$ becomes 
\begin{eqnarray}
\label{deflambda0}
\Lambda_0(T, x, y)  = \int_0^T (1-\tau){L}(X(\tau), Y(\tau))
 d\tau  
\end{eqnarray}
where
\ben
\label{other0}
{L}(X, Y) = {\mathcal J}(1 + X)
- {\mathcal J}(1 + Z) - {\mathcal J}(1+\bar Z)
\een
i.e.
\ben
\nonumber
{L}(X, Y) &=& (1+X)\log (1+X) + Y \arctan \frac{Y}{2+X}\\
\label{other}
&& -\left(1 + \frac{X}{2}\right)\log \left[\left(1 + \frac{X}{2}\right)^2 + \frac{Y^2}{4}\right]\,.
%\label{other}
%&&+ Y \arctan \frac{Y}{2+X}\,.
\een
Looking for the Legendre dual, we see that for $|\eta| < \pi/2$ and $e^\xi < 2 \cos\eta$, the supremum
$${L}^\star(\xi, \eta) = \sup_{X,Y} X\xi + Y\eta - {L}(X,Y)$$
is achieved in $(X,Y)$ satisfying
\ben
\label{eqxieta}
\frac{1+X}{\sqrt{\Big(1+ \frac{X}{2}\Big)^2 + \frac{Y^2}{4}}} = e^\xi\ , \
\frac{Y}{2+X}= \tan \eta
%\frac{\ii}{2}\log\frac{1 + \frac{x - \ii  y}{2}}{1  + \frac{ x + \ii \dot y}{2}} &=& \eta
\een
i.e.
\be
X = \frac{e^\xi - \cos \eta}{\cos \eta - \frac{1}{2}e^\xi}\ , \ Y = \frac{\sin \eta}{\cos \eta - \frac{1}{2}e^\xi}\,.
\ee
Note that under the assumption above, $X$ is admissible, i.e. $X > -1$. One deduces that (cf. (\ref{3.1}))
\ben
\label{3.17}
{L}^\star(\xi, \eta) =  -\xi - \log(2 \cos\eta - e^{\xi}) = H_a(\xi, \eta).
\een
On the other hand, one can check that ${L}^\star(\xi, \eta)$ is infinite if $|\eta| \geq \pi/2$ or $e^\xi \geq 2 \cos\eta$. 

One deduces that there exists a recession function:
\begin{eqnarray*}
(\xi, \eta) \mapsto \lim_{\kappa\rightarrow +\infty} \kappa^{-1}{L^\star}(\kappa \xi, \kappa\eta) = \begin{cases} &-\xi \ \hbox{if} \ \xi<0 \ \hbox{and}\ \eta=0\,,\\
&\infty \ \hbox{otherwise}\,.
\end{cases}
\end{eqnarray*}
%\marginpar{\textbf{to be developped}}
%\[\lim_{\kappa\rightarrow +\infty} \kappa^{-1}{L^\star}(\kappa \xi, \kappa\eta) = -\xi\]
%if $\xi < 0$ and $\eta=0$, and $+\infty$ otherwise.
 This function can be used to obtain the rate function ${\mathcal H}_0(T, \varphi, \psi)$ when 
the measures $\dot \varphi_s$ and $\dot \psi_s$ are not zero. By using the same methods as in \cite{rou}, one deduces
\begin{eqnarray*}
{\mathcal H}_0(T, \varphi, \psi) &=& \sup_{x(.),y(.)} 
\left[\int_0^T \left(x(\tau) d \dot\varphi(\tau) + y(\tau) d \dot \psi(\tau)\right) - \Lambda_0 (T, x(.), y(.))\right]\\%\Lambda(x(\tau), y(\tau))\right)
%=\] 
&=& \sup_{X(.), Y(.)} \int_0^T (1-\tau) \left[X(\tau) d\dot\varphi(\tau) + Y(\tau)d \dot\psi(\tau) -{L}(X(\tau),Y(\tau))\!\ d\tau \right]\\
%{\mathcal H}_0(T, \varphi, \psi)
 &=& \int_0^T (1-\tau) H_a(\dot\varphi_a (\tau), \dot\psi_a (\tau))\!\ d\tau + \int_0^T
(1-\tau) d(-\dot \varphi_s) (\tau)\,,%\,,
\end{eqnarray*}
%\]
where the second equality comes from  (\ref{xtau}) and (\ref{deflambda0}), and where the last equality comes from  \cite{Rocky1} Theorem 5.
%\[{\mathcal H}_0(T, \varphi, \psi) 
%= \sup_{x(.),y(.)} %
%\left[\int_0^T \left(x(\tau) d \dot\varphi(\tau) + y(\tau) d \dot \psi(\tau)\right) - \Lambda_0 (T, x(.), y(.))\right] %\Lambda(x(\tau), y(\tau))\right)
%=\] 
%\[= \sup_{X(.), Y(.)} \int_0^T (1-\tau) \left[X(\tau) d\dot\varphi(\tau) + Y(\tau)d \dot\psi(\tau) -{L}(X(\tau),Y(\tau))\!\ d\tau \right]\,,
%\]
%i.e. (see \cite{Rocky1} Theorem 5)%\marginpar{commentaire}
%\begin{eqnarray*}
%{\mathcal H}_0(T, \varphi, \psi) &=& \int_0^T (1-\tau) H_a(\dot\varphi_a (\tau), \dot\psi_a (\tau))\!\ d\tau + \int_0^T
%(1-\tau) d(-\dot \varphi_s) (\tau)\\
%\end{eqnarray*}

The value of the constant ${\mathcal C}_\dd(T)$ is obtained 
owing to (\ref{lmarg}):
\[{\mathcal C}_\dd(T) = - {\mathcal L}_0(T, 2\Re\!\  \dd, 2\Im\!\ \dd)\,.\]
\medskip

\noindent$\bullet$ \textit{Third step : exponential tightness.}

Exponential tightness is not needed for the second argument since it lives in $[-\pi/2, \pi/2]$.
For the first, we have $|\log x| \leq - \log x + 2 \log 2$ for $x \leq 2$, hence:
\[
\mathbb P_\dd(\sum_{j\leq nT-1} |\log |1-y_j|| \geq na) \leq \mathbb P_\dd( \sum_{j\leq nT-1} -\log (1-y_j) \geq n(a - 2T\log 2)).
\]
Now, for $\theta < 0$ (by Chernov inequality),
\[\mathbb P_\dd( \sum_{j\leq nT-1} -\log (1-y_j) \geq n(a - 2T\log 2))\leq e^{n^2\beta'\theta(a - 2T\log 2)} \mathbb E_\dd\left( |\Phi_{\lfloor nT \rfloor}(1)|^{n^2\beta'\theta}\right) \]
so that, taking logarithm and 
applying (\ref{lmarg}) we get, for $\theta \in (-(1-T) -2\Re\!\ \dd, 0)$
\[\limsup_{n \rightarrow \infty} (\beta' n^2)^{-1} \log \mathbb P_\dd (\sum_{j\leq nT} |\log |1-y_j| \geq na) \leq \theta(a - 2T\log 2) + {\mathcal L}_\dd(T,\theta, 0).\]
It remains to let $a \rightarrow \infty$ to get the exponential tightness.
We remark that when $\dd=0$, the exponential tightness holds only for $T < 1$.
\QED

\begin{rem}The mean trajectory is obtained when $H_a(\dot \varphi, \dot\psi)\equiv 0$ i.e. $\cos \dot \psi = \cosh\dot \varphi$ or $\dot\varphi= \dot\psi = 0$ 
\end{rem}
%\newpage

\subsection{Comment on  Theorem \ref{mastermarg}}% and Theorem \ref{HKOC} } 
%To prove the statement of Theorem \ref{mastermarg}, it is enough to apply the contraction principle to the mapping
%\[(\varphi, \psi) \mapsto (\varphi(T), \psi(T))\,.\]
%which associates the terminal point to a path.
% \QED

Let us study the variational problem (\ref{proj}) issued from the contraction.  Using (\ref{calh0}) and (\ref{3.17}), we see that
the Euler equation is
\begin{eqnarray}\nonumber
\frac{d}{d\tau}\!\ \left((1-\tau)\frac{\partial{L^\star}}{\partial \dot\varphi}\right) = 0\\
\frac{d}{d\tau}\!\ \left((1-\tau)\frac{\partial{L^\star}}{\partial \dot\psi}\right) = 0\,,
\end{eqnarray}
and the optimal path is then given by
\[\dot\varphi(\tau) = \frac{\partial{L}}{\partial X}\left(\frac{\gamma}{1-\tau}, \frac{\rho}{1-\tau}\right)\ , \ \dot\psi(\tau) = \frac{\partial{L}}{\partial Y}\left(\frac{\gamma}{1-\tau}, \frac{\rho}{1-\tau}\right)\]
i.e. 
\ben\label{trajphi}\dot\varphi(\tau) &=&  \log (1 -\tau + \gamma) - \frac{1}{2} \left(\log \Big(1-\tau+ \frac{\gamma}{2}\Big)^2 + \frac{\rho^2}{4}\right),\\
\label{trajpsi}\dot\psi(\tau) &=& \arctan \frac{\rho}{2(1-\tau) + \gamma}.
\een
This path will be admissible if there exist $\gamma$ and $\rho$ such that
%with $\gamma$ and $\rho$ such that
\ben
\label{endpoints}\int_0^T \dot\varphi(\tau) d\tau = \xi \ , \ \int_0^T \dot\psi(\tau) d\tau = \eta\,.\een
When the path is admissible, %\marginpar{what's?},
 we have
\[H_a(\dot\varphi(\tau), \dot\psi(\tau)) = \widehat{\mathcal L} (\dot\varphi(\tau), \dot\psi(\tau)) = \frac{\gamma}{1-\tau}\dot\varphi (\tau) + \frac{\rho}{1-\tau}\dot\psi (\tau)-{L}\left(\frac{\gamma}{1-\tau}, \frac{\rho}{1-\tau}\right)\]
and
\begin{eqnarray*}h_0(T, \xi, \eta) = \gamma \xi + \rho \eta- \int_0^T (1-\tau){L}\left(\frac{\gamma}{1-\tau}, \frac{\rho}{1-\tau}\right)\!\ d\tau.\end{eqnarray*}

%\begin{rem}
\subsection{Comment on Theorem \ref{HKOC}}
\label{HKOC}
In the case $\dd=0$,  $T=1$, \cite{HKOC} proved the LDP by tackling directly the normalized cumulant 
generating function. This gives an incomplete LDP since there is no steepness in $0$. They use a Fourier inversion to take into account the negative side. We see that the function ${\mathcal L}_0(1,s, 0)$ is the  limiting n.c.g.f. %log-Laplace transform
 of $\log |\Phi_n (1)|$, computed in \cite{HKOC} Theorem 3.3. 
Besides, ${\mathcal L}_0(1,0, t)$ is the limiting n.c.g.f. % log-Laplace transform
 of the argument of $\Phi_n(1)$, computed in \cite{HKOC}, Theorem 3.4, 
%\[{\mathcal L}_0(1,0, t) = \frac{t^2}{8}\ln\left(1 + \frac{4}{t^2}\right)\frac{1}{2}\ln\left(1+ \frac{t^2}{4}\right) + t \arctan(t/2)\]
and
\[\lim_{t \rightarrow \pm\infty} \frac{{\mathcal L}_0(1,0, t)}{t} = \pm \frac{\pi}{2}\,.\]
%\end{rem}

\subsection{Proof of Theorem \ref{margmoi}}
In the $\Im\!\ \dd \not= 0$ case, or the case where $\dd = 0$ and $T < 1$, we could also use the
 scheme described in Section \ref{HKOC}. We prefer illustrate the method of contraction, where we will see that the counterpart of the singular contribution is an affine part. To simplify the exposition, we assume the  problem one-dimensional (consider only the first component) and $\dd > 0$.

Since
\be\frac{\partial {\mathcal L}_\dd}{\partial s}(T,s,0) &=& \frac{\partial {\mathcal L}_0}{\partial s}(T,s+2\dd,0)\\
&=& 
{\mathcal J}(1+s+2\dd) - {\mathcal J}(1-T+s+2\dd)\\ && - {\mathcal J}(1 +\frac{s}{2} +\dd) + {\mathcal J}(1 -T +\frac{s}{2}+\dd)\ee
and since ${\mathcal J}(a+s) - {\mathcal J}(b+s) = (a-b) \log s + o(1)$ as $s \rightarrow \infty$, we have
\[\lim_{s \rightarrow \infty} \frac{\partial{\mathcal L}_\dd}{\partial s} (s, 0) = T\log 2\] which corresponds to the endpoint of the interval allowed for
 $\Re \log \Phi_n(T)$.
On the other side, if $\dd > 0$,
\[\lim_{s\downarrow -(1-T)-2 \dd} \frac{\partial{\mathcal L}_d}{\partial s} (T,s, 0) =  {\mathcal J}(T) - 1 - {\mathcal J}\left(\frac{1+T}{2}\right) + {\mathcal J}\left(\frac{1-T}{2}\right) =: \xi_T.\]
There is a problem of non-steepness since it is not $-\infty$.

%Let us consider only the first component to explain the phenomena.
The first coordinate satisfies the LDP with good rate function
\[ \xi \mapsto \inf \{ h_0 (T, \xi, \eta) | \eta\}\]
and from (\ref{proj}) \[ \inf \{ h_0 (T, \xi, \eta) | \eta\}= \inf \{\mathcal H_0(T, \varphi, \psi) | \varphi(T) = \xi\}\]
From the structure of $\mathcal H_0$ and $H_a$, it is clear that this infimum is achieved for $\dot\psi(.) = 0$ i.e. $\rho = 0$.
%\[\]
%Since \[\frac{\partial L}{\partial Y}= \arctan \frac{Y}{2+X}\]
%we see that for $ \varphi$ fixed, $H_a(\dot\varphi, \dot \psi)$ is minimum for $\dot \psi$ such that $\frac{\partial L^\star}{\partial Y} (\dot\varphi, \dot\psi) = 0$, i.e. $\dot\psi = 0$, or in other words $\rho =0$.

In the case of admissible $\varphi$, it  remains to compute $\gamma$. Let us study the mapping  $\gamma \mapsto \varphi(T, \gamma)$ where $\varphi(\cdot, \gamma)$ is given by (\ref{trajphi}) with $\rho =0$. We follow the lines of argument of \cite{rou}  pp. 3216-3218. We have
\ben
\frac{\partial \varphi(T, \gamma)}{\partial \gamma} = \log\frac{1+\gamma}{1-T+\gamma} -\frac{1}{2}\log \frac{1 + \frac{\gamma}{2}}{1-T + \frac{\gamma}{2}} > 0\,. 
\een
The mapping 
$\gamma \mapsto \varphi (T, \gamma)$ is then bijective from $[-(1-T)%-2\dd
,\infty )$ to $[\xi_T, T\log 2)$.

Fixing $\xi \in (-\infty, T\log 2]$, let us look for optimal $\varphi$. Let $\gamma > -(1-T)$ (playing the role of a Lagrange multiplier). By the duality property, we have the inequality 
\[ (1-\tau) {L^\star}\left(\dot \varphi_a(\tau)\right)\geq \gamma \dot \varphi_a(\tau) - (1-\tau) { L}\left(\frac{\gamma}{1-\tau}, 0\right)\,.\]
Using (\ref{calh0}) and (\ref{3.17}) we get, by integration of the above inequality:
%integrating and using (\ref{lncgf}) and the expression of $H_a$, we get
\begin{eqnarray*}
{\mathcal H}_0(T, \varphi, 0) &=& \int_0^T (1-\tau) L^\star(\dot \varphi_a(\tau)) d\tau + \int_0^T (1 - \tau) d \left(-\dot\varphi_s\right)(\tau)\\
&\geq& \gamma \varphi_a (T) - \int_0^T (1-\tau) L\left(\frac{\gamma}{1-\tau}, 0\right) d\tau\\
 &+& \int_0^T  (1 - \tau) d \left(-\dot\varphi_s\right)(\tau).
\end{eqnarray*}
For every $\varphi$ such that $\varphi(T) = \xi$, we then have
\ben
\nonumber
{\mathcal H}_0(T, \varphi, 0)&\geq& \gamma \xi -\int_0^T (1-\tau) L\left(\frac{\gamma}{1-\tau}, 0\right) d\tau -\int_0^T (1-\tau+\gamma) d\dot\varphi_s(\tau)
\\ \label{plug}
 &\geq& \gamma \xi -\int_0^T (1-\tau)L\left(\frac{\gamma}{1-\tau}, 0\right) d\tau.
\een

We can now distinguish three cases:

$\bullet$ If $\xi \in [\xi_T, T\log2)$, we choose the path $v^\xi$ absolutely continuous and
such that $$\dot{v^\xi} (\tau) = \log(1- \tau + \gamma^{\xi}) - \log \left(1 - \tau +
\frac{\gamma{^\xi}}{2} \right),$$
where $\gamma{^\xi}$ is uniquely determined by the condition $v^{\xi} (T) = \xi$.
 This path saturates the infimum and the expression of the action integral is clear.

$\bullet$ If $\xi < \xi_T$, set $\varepsilon = \xi_T-\xi$. Plugging $\gamma= -(1-T)$ in (\ref{plug}) yields for $v$
%$varphi$
 such that $v(T) = \xi$:
 %$\varphi(T) =\xi$
\[{\mathcal H}_0(T, v, 0) \geq -(1-T) \xi -\int_0^T (1-\tau)L\left(\frac{-(1-T)}{1-\tau}, 0\right)= (1-T)\varepsilon +h_0(T, \xi_T, 0)\]
and this lower bound is achieved for the measure
$\widetilde v = v^{\xi_T}(\tau) d\tau -\varepsilon \delta_T$, since
\[{\mathcal H}_a (v^{\xi_T}) = h_0  (T, \xi_T, 0) \ , \ \int_0^T (1-\tau)\varepsilon\!\ d\delta_T (\tau) = (1-T)\varepsilon\,.\]

$\bullet$ If $\xi = T \log 2$, make $\xi = T \log 2$ in (\ref{plug}). We get, for 
all $\gamma > -(1-T)$,  
\ben\label{lastineg}h_0 (T, T\log 2,0) \geq  \gamma T\log 2 - \int_0^T (1-\tau) L\left(\frac{\gamma}{1-\tau}, 0\right) d\tau. \een
When $\gamma \rightarrow \infty$, the integral is $F(1+\gamma) - F(1-T+\gamma) + F(1) 
-  F(1 - T) -2 F(1+\gamma/2) + 2 F(1-T+\gamma/2)$ which 
tends to $-\infty$ as $-T \log \gamma$, so that finally, the RHS (\ref{lastineg}) tends to $\infty$ and we conclude $h_0 (T, T\log 2,0)
= \infty$.
\section{Connection with the spectral method}

It can be interesting to connect the results of the previous section to the results obtained
by looking directly at the empirical spectral distribution of the ensembles which are considered. 
This point of view is also discussed in \cite{HKOC}.

The LDP for the empirical spectral distribution of the unitary ensemble is given in \cite{HiaiIHP}. The rate function is the Voiculescu's logarithmic entropy:
\ben \label{HP2} I(\mu) = - \Sigma_{\TT}(\mu) :=
- \int\!\!\int_{\mathbb T \times \mathbb T} \log|z-z'|\ d\mu(z) d\mu(z'). \een
The circular Jacobi unitary ensemble yields also a LDP given in \cite{apa}. The rate function is
\ben
\label{sigmaq}I_\dd(\mu) = %- \Sigma_\dd(\mu)
  - \Sigma_{\TT}(\mu) +\int_\mathbb T  Q_\dd(z) d\mu(z) + B(\dd)\,,\een
where
\begin{equation}
Q_\dd(e^{\ii \theta}):= - 2(\Re\!\ \dd) \log \big(2 \sin \frac{\theta}{2}\big) - (\Im\!\ \dd) 
(\theta - \pi) \ \ (\theta \in (0, 2\pi))\,.
\end{equation}
and
\begin{eqnarray}
\nonumber
B(\dd) &=& \int_0^1 \left[(x+ 2\Re\!\ \dd)\log (x + 2\Re\!\ \dd)-2\Re\!\ [(x+  \dd)\log (x +  \dd)]\right] dx\\
&+&\label{vraieconstante} \int_0^1 x \log  \dd x\,.\end{eqnarray}
(Notice that there was a mistake in \cite{apa}, fixed in the arXiv version). 
With our notation, it yields:
\begin{equation*}
%\label{autre}
B(\dd) = F(1 + 2\Re\!\ \dd) - F(2\Re\!\ \dd) - 2 \Re\!\ F(1 + \dd) +2 \Re\!\ F(\dd) + F(1)\,.
\end{equation*}
%\[B(\dd) = \int_0^1 x \log \frac{x(x + \Re\!\  \dd)}{|x+\dd|^2}\!\ dx\,.\]

If the mapping $\mu \in {\mathcal M}_1 ({\TT}) \mapsto \int \log (1-z)\!\ d\mu(z)$ were continuous, we would have by contraction:
\ben
\label{vraie}
h_\dd(1,\xi, \eta) = \inf \{I_\dd (\mu) %- \Sigma_{\dd}(\mu)
 \ |\  \mu \in {\mathcal M}_1(\TT) : \ \int_\mathbb T \log(1-z)\!\   d\mu(z) = \xi + \ii \eta\}\,.
\een
We conjecture that this formula holds nevertheless, and we prove it in the one-dimensional case. 

\begin{prop}
\label{conject}
With the notation of Section 3, we have, for $\dd > 0$ and for every $\xi \in [0, \log 2)$
\ben
\label{conject1}
h_\dd(1,\xi, 0) = \inf \{%-\Sigma_{\TT}
 I_\dd(\mu) \ | \ \int_{\TT} \log|1-z|\!\   d\mu(z) = \xi \}\,.
\een
\end{prop}

We use the following interesting result.
\begin{prop}
\label{minr}
We have
\ben
\inf \{- \Sigma_\TT (\mu) \ |\  \mu \in {\mathcal M}_1(\TT) :  \ \int_\TT \log |1-z|\!\   d\mu(z) = \xi\} = - \Sigma_\TT (\mu_a)
\een
where 
\begin{itemize}
\item the measure $\mu_a \in {\mathcal M}_1(\TT)$ is defined by
\begin{equation}
\label{limmeas}
d\mu_a(z) = (1+a) \frac{\sqrt{\sin^2\big(\theta/2\big) -
\sin^2(\theta_a/2)}}{2\pi \!\ \sin(\theta/2)} \mathds{1}_{(\theta_a, 2\pi-\theta_a)}(\theta)\!\   d\theta
\end{equation}
where $z= e^{\ii \theta}, \theta \in [0, 2\pi]$ and 
$
\theta_a \in (0, \pi)$ is such that $\sin \theta_a/2 = \frac{a}
{ 1 + a}$.
%\item
%\[d\nu_a(z) = \frac{\Big(1 + \sqrt{1+a^2}\Big)}{a\pi}\!\ \sqrt{1 - \frac{1}{a^2}\!\ \left|\frac{z+1}{z-1}\right|^2}\!\   {\mathbb I}_{z \in {\mathcal %C}_a} dz\]
%\item $dz$ is the Haar measure on $\TT$,
\item $a$ is the unique solution of
\[\int_\TT \log|1-z|\!\  d\mu_a(z) = \xi.\]
%or, equivalently
%\[{\mathcal J}(1+2a) - {\mathcal J}(2a) -{\mathcal J}(1+a) + {\mathcal J}(a) = \xi\]
%\item 
%${\mathcal C}_a = \{ z \!\ ; z  = e^{\ii \theta} \ , |\tan(\theta/2)| \geq a^{-1} \}$
\end{itemize}
Moreover,
\[\int_\TT \operatorname{arg}(1-z)\!\   d\mu_a(z) = 0\,.\]
\end{prop}

The two propositions mean  that to get the non standard mean value $\xi$ for the logarithm of the determinant, the most probable way is to force the random operator to get an empirical spectral distribution close to $\mu_a$.
\medskip

\noindent\textsl{Proof of Proposition \ref{conject} given Proposition \ref{minr}.}
Let us first suppose that the Proposition \ref{conject} is true for $\dd = 0$. One deduces, from this assumption and the definition (\ref{sigmaq}) of $I_\dd$:
\begin{eqnarray*}
&&\inf \{%-\Sigma_{\TT}
 I_\dd(\mu) \ |\  \mu \in {\mathcal M}_1(\TT) : \ \int \log|1-z| d\mu(z) = \xi \}\\
%\\  &=& \inf \{-\Sigma_{\TT} (\mu)
% + \int Q_{\dd} (\theta) d \mu(\theta)  + B(\dd) 
 % \ | \ \int \log|1-z| d\mu(z) = \xi \}
&&=      \inf \{-\Sigma_{\TT} (\mu) \ |\  \mu \in {\mathcal M}_1(\TT) : \ \int \log|1-z| d\mu(z) = \xi \}
  + B(\dd) - 2d \xi   \\&& =h_0(1, \xi, 0) + B(\dd) -  2\dd \xi = h_{\dd} (1, \xi, 0),
\end{eqnarray*}
the last equality coming from \eqref{0d} and the fact that $B(\dd) = - C_{\dd} (1)$. 
Hence, it is sufficient to prove Proposition \ref{conject} for $\dd = 0$. 

The RHS of (\ref{conject1}) is $- \Sigma_\TT (\mu_a)$ and it
 remains to prove that it fits with the value
\[h_0(1,\xi, 0) = \gamma\xi -
F(1+\gamma) + F(\gamma) - F(1) +2 F\left(1 + \frac{\gamma}{2}\right) -2F\left(\frac{\gamma}{2}\right)\]
where in view of Theorem \ref{margmoi} (2), $\gamma$ is the solution of (\ref{implicit}) i.e.
\ben
{\mathcal J}(1+\gamma) - {\mathcal J}(\gamma) - {\mathcal J}\left(1 + \frac{\gamma}{2}\right) + {\mathcal J}\left( \frac{\gamma}{2}\right)= \xi\,.
\een
Applying \cite{apa} formula (5.25) in the special case $\dd = a$ (real and positive), we see that,
since the corresponding rate function vanishes for the limiting empirical measure $\mu_a$,
\ben
\Sigma_{\TT}(\mu_a) = \int_\TT Q_a(\zeta)\!\   d\mu_a (\zeta) + B(a)\,,
\een
where
$Q_a(\zeta) = -2 a \log |1- \zeta|$.
Let us now compute:
\begin{equation}
\label{mys}
I:= \frac{2\pi}{1+a}\int_\TT \log |1-\zeta|\ d\mu_a(\zeta)\,. 
\end{equation}
We have, by definition of $\mu_a$ in (\ref{limmeas}),
\begin{eqnarray}
\nonumber
 I &=&  \int_{\theta_a}^{2 \pi - \theta_a}  \log (2 \sin (\theta/2)) \frac{\sqrt{\sin^2(\theta/2) - \sin^2(\theta_a/2)}}
{\sin(\theta/2)} \, d\theta\\ 
\label{mys2}
&=&
 2  \int_{\theta_a}^{\pi}  \log (2 \sin (\theta/2)) \frac{\sqrt{\sin^2(\theta/2) - \sin^2(\theta_a/2)}}
{\sin(\theta/2)} \, d\theta\,,
\end{eqnarray}
for $\theta_a \in [0, \pi]$ and $\sin(\theta_a/2) = a/(1+a)$.
The first change of variables $u = \sin(\theta/2)$  gives 
$$d\theta = \frac{2 du}{\cos(\theta/2)} = \frac{2 du}{\sqrt{1 - u^2}},$$
and then
$$ I = 4 \int_{a/(1+a)}^1 \log (2u) \sqrt{\frac{u^2 - [a/(1+a)]^2}{1 - u^2}} \frac{du}{u}.$$
The second change of variable
$$v = \sqrt{\frac{u^2 - [a/(1+a)]^2}{1 - u^2}} \Longleftrightarrow
 u = \sqrt{\frac{v^2  + [a/(1+a)]^2}{1 + v^2}},$$
gives
$$\log (2u) = \log 2 + \frac{1}{2} \log \left(v^2  + [a/(1+a)]^2 \right)
- \frac{1}{2} \log (v^2  + 1 ),$$
$$\frac{du}{u} = \left[ \frac{v}{v^2  + [a/(1+a)]^2} - \frac{v}{v^2 + 1}  \right]\, dv.$$
Hence, 
$$ I = 4 \int_0^{\infty} \left[ \log 2 + \frac{1}{2} \log \left(v^2  + [a/(1+a)]^2 \right)
- \frac{1}{2} \log (v^2  + 1 ) \right] \cdot  \left[ \frac{v^2}{v^2  + [a/(1+a)]^2} - \frac{v^2}{v^2 + 1}  \right]
\, dv$$
\begin{equation} = 2 \int_{-\infty}^{\infty} \left[ \log 2 + \frac{1}{2} \log \left(v^2  + [a/(1+a)]^2 \right)
- \frac{1}{2} \log (v^2  + 1 ) \right] \cdot  \left[ \frac{1}{v^2 + 1} -
 \frac{[a/(1+a)]^2}{v^2  + [a/(1+a)]^2} \right]
\, dv. \label{expressionI}
\end{equation}
Now, for $\alpha, \beta > 0$, 
\begin{equation}\int_{-\infty}^{\infty} \frac{\alpha^2}{v^2 + \alpha^2} dv = \alpha \pi, \label{integral1}
\end{equation}
\begin{equation}
\int_{-\infty}^{\infty} \frac{\alpha^2 \, \log (v^2 + \beta^2)}{v^2 + \alpha^2} dv = 2 \alpha \pi \log(\alpha + \beta).
\label{integral2}
\end{equation}
The first equality \eqref{integral1} is elementary. To get  \eqref{integral2}, one can 
differentiate the LHS with respect to $\beta$, and one obtains, for $\alpha \neq \beta$,
\begin{eqnarray*}
 \int_{-\infty}^{\infty} \frac{2 \alpha^2 \, \beta}{(v^2 + \alpha^2)(v^2 + \beta^2)} dv
=  \frac{2  \alpha^2 \, \beta}{\beta^2 - \alpha^2} \left[ \int_{-\infty}^{\infty} \frac{dv}{v^2 + \alpha^2}
-  \int_{-\infty}^{\infty} \frac{dv}{v^2 + \beta^2} \right]
%\\ &=&  \frac{2  \alpha^2 \, \beta}{\beta^2 - \alpha^2} \left( \frac{\pi}{\alpha} - \frac{\pi}{\beta} \right)
%\\ & =  \frac{2 \pi  \alpha^2 \, \beta}{\beta^2 - \alpha^2} \, \frac{\beta - \alpha}{\alpha \beta}
 = \frac{2 \pi \alpha}{\alpha + \beta}.
\end{eqnarray*}
By continuity, the last equality remains true for $\alpha = \beta$, and reversing the differentiation 
gives
$$\int_{-\infty}^{\infty} \frac{\alpha^2 \, \log (v^2 + \beta^2)}{v^2 + \alpha^2} dv = 2 \alpha \pi \log(\alpha + \beta)
+ C(\alpha).$$
Now, to determinate $C(\alpha)$ let us write
$$\int_{-\infty}^{\infty} \frac{\alpha^2 \, \log (v^2 + \beta^2)}{v^2 + \alpha^2} dv 
= \alpha \log(\beta^2) \int_{-\infty}^{\infty} \frac{\alpha}{v^2 + \alpha^2} dv 
+ \alpha \int_{-\infty}^{\infty} \frac{\alpha  \log (1 + v^2/\beta^2)}{v^2 + \alpha^2} dv\,.  $$
The first integral is equal to $2 \alpha \pi \log(\beta)$ and the second tends to zero for $\alpha$ fixed and $\beta \rightarrow \infty$
by dominated convergence. This gives $C(\alpha) = 0$ and then \eqref{integral2}.
Expanding \eqref{expressionI}, then using \eqref{integral1} and \eqref{integral2} gives:
\begin{eqnarray*}
I = % 2 \pi \log(2) - 2 \pi [a/(1+a)] \log(2) + 2 \pi \log [(1 + 2a)/(1+a)] \\ &-& 
%2 \pi [a/(1+a)] \log [2a/(1+a)] - 2 \pi \log(2) + 2 \pi [a/(1+a)] \log [(1 + 2a)/(1+a)] 
%\\ &=& \frac{2 \pi}{1+a} \left[ - a \log (2) + (1+a) \log (1+2a) - (1+a) \log (1+a)\right]  \\ &-& a \log (2a) 
%+ a\log(1+a) + a \log(1+2a) - a \log(1+a) \right]
%\\ &=& \frac{2 \pi}{1+a} \left[ - a \log (2) + (1+2a) \log (1+2a) - (1+a) \log (1+a)  - a \log (2a) \right]
 \frac{2 \pi}{1+a} \left[ (1+2a) \log (1+2a) - (1+a) \log (1+a)  - 2a \log (2a) + a \log(a) \right].
\end{eqnarray*}
Coming back to the definition of $I$, we obtain:
\begin{equation*}
%\label{idef}
\int_\TT \log |1-\zeta|\ d\mu_{a}(\zeta) = \mathcal I(1+2a) -\mathcal I(1+a) -\mathcal I(2a) + \mathcal I(a) = \mathcal E_a(1)\,.
\end{equation*}
In particular, 
$$\int_\TT \log |1-\zeta|\ d\mu_{\gamma/2}(\zeta) = \xi,$$
and the value of $a$ involved in Proposition \ref{minr} is equal to 
$\gamma/2$. It remains to check that 
$- \Sigma_\TT (\mu_{\gamma/2}) = h_0(1,\xi, 0)$. Indeed, 
\be- \Sigma_{\TT}(\mu_{\gamma/2}) &=& - \int Q_{\gamma/2}(\zeta) d\mu_{\gamma/2} (\zeta) - B(\gamma/2)
=  \gamma \xi - B (\gamma/2) \\ &=& \gamma \xi - F(1+\gamma) + F(\gamma)
+ 2F(1 + \gamma/2) - 2 F(\gamma/2) - F(1)\,,\ee
which proves Proposition \ref{conject}. 
% Finally, choosing $\dd = \gamma/2$ yields
%\[\int \log |1-\zeta| d\mu_\dd (\zeta) = \xi\,,\]
%so that
%$\mu_\dd = \nu_a$ with $a = 2(1+\dd)/\dd$ and
%\be -\Sigma(\nu_a) = -\Sigma(\mu_\dd) = \gamma\xi -B(\dd) h_0(1,\xi, 0)\,.
%\ee
\QED
\begin{rem}
%We can try a direct analytical method, but we prefer a probabilistic point of view. 
The determination of $I$ could also  be viewed as a consequence of the weak convergence of the empirical spectral measure under CJ$_{n a \beta'}$. Of course the mapping $\mu \mapsto \int_\mathbb T \log |1-\zeta|\!\   d\mu(\zeta)$ is not continuous, but since the extremal eigenvalues converge to the extreme points of the support of $\mu_a$ (see \cite{Borot}), we have
\begin{equation*}
\label{idefa}
\int_\TT \log |1-\zeta|\!\ d\mu_a(\zeta) =
\lim_{n \rightarrow \infty} \frac{1}{n} \log |\Phi_n (1)| = {\mathcal E}_a (1)%\\ \nonumber &=& {\mathcal J}(1+2a) - {\mathcal J}(2a) -{\mathcal J}(1+a) + {\mathcal J}(a)
\end{equation*}
where the  last equality comes from (\ref{2.7})  in Theorem \ref{2.1}.
\end{rem}
\bigskip

\noindent \textsl{Proof of Proposition \ref{minr}} A possible way consists in studying the maximization of the functional
\begin{eqnarray*}
\Sigma_{\TT}(\mu) +r\int_\TT \log|1-z|\ d\mu(z)\,,
\end{eqnarray*}
where $r$ is some Lagrange multiplier.

 It is convenient to push-forward this problem to $\mathbb R$ via the Cayley transformation
\[ z = e^{\ii \theta} = \frac{ \lambda +\ii}{\lambda -\ii}\]
so that, for $z, z' \neq 1$, \ben\nonumber\log|1-z| &=& -\frac{1}{2}\log(1 + \lambda^2) + \log 2 
\\\nonumber \log|z-z'| &=& \log|\lambda - \lambda'| - \frac{1}{2}\log (1+ \lambda^2) - \frac{1}{2}\log (1+ \lambda'^2) + \log 2 \,.\een Hence, one needs to minimize:
\[\Sigma_{\mathbb R}(\nu) + 2\int_{\mathbb R} Q(\lambda)\!\ d\nu(\lambda),\]
over $\nu \in \mathcal M_1(\mathbb R)$, where
\[ \Sigma_{\mathbb R}(\nu) = \int\!\!\int_{\mathbb R \times \mathbb R} \log \left( \frac{1}{|\lambda - \lambda'|} \right)
 d \nu(\lambda) d \nu(\lambda')
\]
and
\[2 Q(x) = \big(1+\frac{r}{2}\big)  \log (1+ x^2). 
\]
This problem can be connected to the research of the equilibrium measure in the Cauchy ensemble, i.e. for the Coulomb gas with probability distribution given
by formula (3.124) in  Forrester \cite{Forrester}.
According to Saff and Totik \cite{Saff}, this potential $Q$ is admissible (p.27) as soon as $r> 0$. In that case the minimizer is unique (Theorem 1.3), its support is compact.
In the book, an explicit method is presented for a class  of potentials satisfying some conditions,
in particular if $Q$ is even, differentiable and such that $xQ'(x)$ is positive and increasing in $(0, \infty)$. These properties are satisfied here since: 
\[xQ'(x) = \big(1+\frac{r}{2}\big)\frac{x^2}{1+x^2}\,.\]
Then (\cite{Saff} Corollary 1.12 p. 203) the support is $S = [-b, b]$ where $b$ is solution of
\ben\label{endpoint}\frac{2}{\pi}\int_0^1 \frac{btQ'(bt)}{\sqrt{1-t^2}}\!\ dt = 1\een
i.e.
\[\frac{2+r}{\pi}\int_0^1 \frac{b^2t^2}{(1+ b^2t^2)\sqrt{1-t^2}}\!\ dt = 1 
\]
or
\ben
\label{aux}\int_0^1 \frac{dt}{(1+ b^2t^2)\sqrt{1-t^2}}= 
\frac{\pi r}{2(2+r)} \een
which gives (make $t= (1 + \tau^2)^{-1/2}$)
\[
\int_0^\infty \frac{d\tau}{1+b^2 + \tau^2}= 
\frac{\pi r}{2(2+r)} \]
i.e. \[b = \frac{2\sqrt{1+ r}}{r}\,.\]
To find the extremal measure, we apply Theorem IV.3.1 in \cite{Saff},
which contains the following result:
\begin{thm} [Lubinsky-Saff] \label{LS2011}
Let $f$ be  a differentiable even function on $[-1, +1]$,  such that $sf'(s)$ is increasing for $ s \in(0,1)$ and for some $1 < p < 2$, $f'(s)/\sqrt{1-s^2} \in L^p[-1, 1]$. Then the integral equation
\[\int_{-1}^1 \log\frac{1}{|x-t|} g(t) dt = -f(x) +C_f, \ \ \ x \in (-1, 1)\,,\]
(where $C_f$ is some constant) 
has a solution of the form
\ben
g(t) = L[f'](t) + \frac{B_f}{\pi\sqrt{1-t^2}}
\een
where
\ben
L[f'](t) &=& \frac{2}{\pi^2}\int_0^1 \frac{\sqrt{1-t^2}\!\ (sf'(s) -tf'(t))}{\sqrt{1-s^2}\!\ (s^2-t^2)}\!\ ds,\\
\nonumber
B_f &=& 1 - \frac{1}{\pi}\int_{-1}^1 \frac{sf'(s)}{\sqrt{1-s^2}}\!\ ds\,.
\een
\end{thm}
 Here, we look for the measure $\mu$ with support $[-b,b]$ such that for $z \in (-b, b)$,
\ben
\label{nonscaled}
\int_{-b}^b \log \frac{1}{|z-t|} d\mu(t) = - Q(z) + C.\een
Theorem \ref{LS2011} is set for  $b=1$, so we need a scaling. 
Equation (\ref{nonscaled}) becomes, with $d\mu(x) = g_b(x) dx$,
\ben
\label{scaled}
\int_{-1}^{1} \log \frac{1}{|z-t|} bg_b(bt) dt  = - Q(bz) + C'.\een
Theorem \ref{LS2011} can now be applied by taking $f(s) = Q(bs)$, and one has
$b g_b(bt) = g(t)$, $B_f = 0$ and
\[\frac{sf'(s) - tf'(t)}{s^2-t^2}= b^2 \Big(1 + \frac{r}{2}\Big) \frac{1}{(1 + b^2t^2)(1+b^2s^2)}\,.\]
Owing to (\ref{aux}), we get
\[bg_b(bt) = \frac{\Big(1 + \sqrt{1+b^2}\Big)}{\pi}\!\ \frac{\sqrt{1-t^2}}{(1 + b^2t^2)}\!\  {\mathbb I}_{[-1, 1]}(t)\]
i.e.
\[g_b (x) = \frac{\Big(1 + \sqrt{1+b^2}\Big)}{b\pi}\!\ \frac{\sqrt{1-x^2b^{-2}}}{(1 + x^2)}\!\  {\mathbb I}_{[-b, b]}(x).\]
Now, it is enough to apply Theorem I.3.3 of \cite{Saff} as follows. The support of the equilibrium measure is $S = [-b, b]$, the measure $\mu$ is supported by $S$, has a finite logarithmic energy and  $\int_{-b}^b \log \frac{1}{|z-t|} d\mu(t) + Q(z)$ is constant for $z \in S$, so $\mu$ is the equilibrium measure. 
To find $\mu_a$, it is enough to  carry $\mu$  back on the circle by the inverse  transformation
\[\lambda= \frac{1}{i} \frac{1+z}{1-z}\,.\]
 \QED
  
\section{Appendix}
\subsection{Some properties of $\ell = \log \Gamma$ and $\Psi = \ell'$}
\label{appendix} From the Binet formula (Abramowitz and  Stegun
\cite{astig} or Erd\'elyi et al. \cite{bateman} p.21), we have for $\Re\!\ x > 0$ 
\ben
\label{bin1}
\ell (x) &=& (x-\frac{1}{2})\log x -x +1 + \int_0^\infty f(s)[e^{-sx} - e^{-s}]\ \! ds\\
\label{bin2}
&=& (x-\frac{1}{2})\log x -x + \frac{1}{2} \log(2\pi) + \int_0^\infty f(s)e^{-sx}\ \! ds\,,
\een
where the function $f$ is defined by
\ben
\label{propf}
 f(s) = \left[\frac{1}{2}-\frac{1}{s}+ \frac{1}{e^s -1}\right]\frac{1}{s} =
 2\sum_{k=1}^\infty \frac{1}{s^2 + 4\pi^2 k^2}\,,
\een and satisfies for every $s \geq 0$: \ben \label{proprf} 0 < f(s)
\leq f(0)= 1/12 \ , \  \ 0 < \left(sf(s)+ \frac{1}{2}\right) < 1\,.
\een By differentiation (recall that $\Psi$ is the Digamma function $\Gamma'/\Gamma$),
\ben
\label{bin3}
\log x - \Psi (x) = \frac{1}{2x} +\int_0 ^\infty s f(s) e^{-sx}\, ds =
 \int_0 ^\infty e^{-sx} \left(sf(s)+ \frac{1}{2}\right)
 ds\,.
\een
Moreover, since $\Psi(x+1) = \frac{1}{x} + \Psi(x)$, we have a variation of (\ref{bin3}):
\begin{eqnarray}
\label{bin4}
\log x - \Psi (x+1) =- \frac{1}{2x} +\int_0 ^\infty s f(s) e^{-sx}\, ds\\
\nonumber
 = \int_0 ^\infty e^{-sx} \left(sf(s)- \frac{1}{2}\right)
 ds\,.
\end{eqnarray}
As easy consequences, we have, for every $x > 0$,
\ben
\label{supx}
0 &<& x\left(\log x - \Psi (x)\right) \leq 1\,,\\
\label{supx2} 0 &<& \log x - \Psi(x) - \frac{1}{2x} \leq
\frac{1}{12x^2}\,. \een Differentiating again we see that for $q\geq
1$, $\Re\!\ x > 0$, 
\begin{eqnarray*} 
%\label{polygamma}
 \Psi^{(q)}(x) = (-1)^{q-1} (q-1)! x ^{-q} +
(-1)^{q-1} \int_0^\infty e^{-sx} s^q  \left(sf(s)+
\frac{1}{2}\right)\ ds \end{eqnarray*}
 and then \ben \label{restepsi}
 |\Psi^{(q)}(x) - (-1)^{q-1} (q-1)! x ^{-q}| \leq (\Re\!\ x)^{-q-1} q!\,.
 \een
Another useful formula is
\ben
\label{Psisum}
\Psi(z+1) = -\gamma - \sum_{k=1}^\infty \left(\frac{1}{k+z}-\frac{1}{k}\right)\,,
\een
fir $z + 1 \notin \mathbb{R}_-$. 
\subsection{The density $g_r^{(\delta)}$ and some moments related to it}
Recall that for $r >0$ and $\delta$ such that $r + 2\Re\!\ \delta +1 >0$, the density $g_{r}^{(\delta)}$ on the unit disc $\mathbb D$  is given by
\[g_{r}^{(\delta)}(z) = c_{r, \delta}  (1 -|z|^2)^{r -1} (1-z)^{\bar\delta} (1- \bar z)^\delta \]
where $ c_{r, \delta}$ is the normalization constant. The following lemma is the key to
 compute $ c_{r, \delta}$ and the moments of $g_r^{(\delta)}$.
\begin{lem}\label{lemhypergeo}
Let $s,t,\ell$ be complex numbers such that: $\Re\!\ \ell$, $\Re (s + \ell +1)$,
 $\Re( t +\ell + 1 )$ and $\Re( s + t +\ell + 1 )$ are strictly positive. Then, the following identity holds:
\ben
\label{laphua}
\int_{\mathbb D} (1 -|z|^2)^{\ell -1} (1-z)^s (1- \bar z)^t \dd^2z = 
\pi  \Gamma\left[\begin{matrix}\let \over/\ell\ ,\ \ell+1+s+t\\ \let\over / \ell +1+s\ , \ \ell+1+t \end{matrix}\right]\,,
%\frac{\pi \Gamma(\ell)\Gamma(\ell +1 + s+t)}{\Gamma(\ell+1+s) \Gamma(\ell + 1+t)}\,. 
\een
where for the sake of simplicity we use the polygamma symbol
\[\Gamma\left[\begin{matrix}\let \over/ a, b, \cdots\\ \let \over /c,d, \cdots\end{matrix}\right] := \frac{\Gamma(a)\Gamma(b)\cdots}{\Gamma(c)\Gamma(d)\cdots}.\] 
\end{lem}
A proof of this result is given in \cite{apa}. 
A first consequence is that
\ben
\label{firstc}
c_{r, \delta}= \pi^{-1}  \Gamma\left[\begin{matrix}\let \over / r +1+\delta\ , \ r+1+\bar\delta\\ \let\over/r\ ,\ r+1+\delta+ \bar\delta\\ \end{matrix}\right]\,.
\een
A second consequence is that if $\gamma$ has the density $g_{r}^{(\delta)}$, then we have
\ben
\nonumber
%\label{mfy}
\mathbb E (1- \gamma)^a (1- \bar \gamma)^b =
 \Gamma\left[\begin{matrix}\let \over/ r+1+\delta+\bar\delta+a+b\ ,\  r + 1 + \bar{\delta}\ ,\  r+1+\delta\\ \let \over  
/  r+1+\delta+\bar\delta\ ,\  r+1+\bar\delta+ a \ ,\ r+1+\delta+b \end{matrix}\right]
\\ \label{mfy}
\een
as soon as all the real parts of the arguments of the gamma functions are strictly positive. 

Let us notice that for $r=0$ the RHS of (\ref{mfy}) is the Mellin-Fourier transform of $1-\gamma$ when 
$\gamma \in \TT$ is
 distributed according to $\lambda^{(\delta)}$.

In this paper, we need the following computations, in order to deduce the moments of $\log(1-\gamma)$. 
The quantities involved below are all well-defined as soon as $s > s_0$, where $s_0$ is some strictly 
negative quantity depending on $r$ and $\delta$, and in particular, for $(s,t)$ in the neighborhood of 
$(0,0)$, one can write: 
\begin{eqnarray}
\nonumber 
\Lambda(s,t)& :=& \log \mathbb E  \exp\left(2s \Re \log (1-\gamma) +2 t \Im \log(1-\gamma)\right) =\\
\nonumber 
&=& \log \mathbb E \exp \left(\Re\!\ (2(s-\ii t )\log (1-\gamma)\right)=\\
\nonumber 
 &=&\log \mathbb E (1-\gamma)^{s-\ii t} (1-\bar \gamma)^{s+\ii t}=\\
\label{core}
&=&  \ell \big(r +1+\delta +\overline\delta+ 2s\big)- \ell\big(r+1+ \delta +\overline\delta\big)\\
\nonumber 
&& - \ell\left(r +1+ \overline\delta+s-\ii t\Big) - \ell\Big(r +1+ \delta+ s+\ii t\right) \\
%\label{core}
\nonumber &&+\ell \big(r +1 +\overline\delta\big)+ \ell\big(r +1 +\delta\big).
\end{eqnarray}
%\proof
To compute  moments we need differentiation. First we have:
\begin{eqnarray}
\nonumber
\frac{\partial}{\partial s}\Lambda(s, t) &=& 2 \Psi \left(r +1+\delta +\overline\delta+2s\right)\\
%\nonumber
\label{once}
&& -\Psi \left(r +1+ \delta+ s+ \ii t\right)-\Psi\left(r +1+ \overline\delta+s- \ii t\right)\\
\nonumber
\frac{\partial}{\partial t}\Lambda(s, t) &=& \ii\Psi \Big(r +1+\bar \delta+ s -\ii t\Big) - \ii \Psi \Big(r +1+\delta +s + \ii t\Big)\,.
%\\
%\label{once}
\end{eqnarray}

The first  moment is then:
\be\mathbb E\!\  \Re \log (1-\gamma) &=& \Psi(r+1 + \delta + \bar\delta) - \frac{1}{2}\Psi(r+1+\delta) -\frac{1}{2}\Psi(r+1+\bar\delta)\\
\mathbb E\!\  \Im \log(1-\gamma) &=&  \frac{1}{2\ii}\Psi(r+1+\delta) -\frac{1}{2\ii}\Psi(r+1+\bar\delta)\ee
or 
\ben
\label{m1c}
\mathbb E \log (1-\gamma) =  \Psi(r+1 + \delta + \bar\delta) - \Psi(r+1+\bar\delta)\,.
\een

Differentiating again (\ref{once}) we get
\begin{eqnarray}
\nonumber
\frac{\partial^2}{\partial s^2} \Lambda(s, t)&=&4 \Psi' \big(r +1+\delta +\overline\delta+2s\big)\\
\nonumber
&&- \Psi'\Big(r +1+ \delta+ s + \ii t\Big)-\Psi'\Big(r +1+ \overline\delta+s- \ii t\Big)\\
 \label{twice}\frac{\partial^2}{\partial t^2}\Lambda(s, t) &=& \Psi' \Big(r +1+\bar \delta+s-\ii t\Big) +\Psi' \Big(r +1+\delta +s +\ii t\Big)\\
\nonumber\frac{\partial^2}{\partial s\partial t} \Lambda(s, t)&=&-\ii\Psi'\Big(r+1+\delta+s + \ii t\Big) +\ii\Psi'\Big(r+1+\bar\delta+s-\ii t\Big)
\end{eqnarray}
and the second moments are
\be
\hbox{Var}\!\ \Re \log (1-\gamma) = \Psi'
 \big(r +1+\delta +\overline\delta\big) -\frac{1}{4}\Psi'\big(r +1+ \delta\big)-\frac{1}{4}\Psi'\big(r +1+ \overline\delta\big)
\ee
\ben
\label{covtheo}
 \hbox{Var}\!\ \Im \log(1-\gamma) = \frac{1}{4}\Psi' \big(r +1+\delta\big) + \frac{1}{4}\Psi' \big(r +1+\overline\delta\big)\een
\be
\hbox{Cov}\!\  (\Re \log (1-\gamma), \Im \log(1-\gamma)) = \frac{1}{4\ii}\Psi'\big(r+1+\delta\big) - \frac{1}{4\ii}\Psi'\big(r+1+\bar\delta\big).
\ee

\subsection{Complex logarithm and characteristic polynomial}

Let $E_k$ be the set of the complex $k \times k$ matrices
with no eigenvalue on the interval $[1, \infty)$. For $V \in E_k$, let us define 
$$\log \operatorname{det} (I_k - V) := \sum_{j=1}^k \log (1 - \lambda_j),$$
where the $\lambda_j$'s are the roots, counted with multiplicity, of the polynomial $z \mapsto \det(z I_k -V)$, and where 
in the right-hand side, one considers the principal branch of the logarithm.
 This definition is meaningful, since by assumption, $1 - \lambda_j \notin \mathbb{R}_-$ for all $j \in \{1, \dots, k\}$. 
 By the continuity of the set of roots of a polynomial with respect to its coefficients, the set $E_k$ is open 
and the function $V \mapsto \log \operatorname{det} (I_k - V)$ defined just above is continuous on $E_k$. In fact, since $E_k$ is connected (this is 
easily checked by tridiagonalizing  the matrices), this is the 
unique way to define the logarithm of $\operatorname{det} (I_k - V)$ as a continuous function of $V \in E_k$ if we assume that it 
should take the value zero at $V =0$. 

Now, with the notation of the beginning of the paper, the matrix $G_k(U_n)$ is a submatrix 
of the unitary matrix $U_n$, and all its eigenvalues have modulus bounded by $1$. 
If we assume $ \gamma_0, \dots, \gamma_{n-1} \neq 1$ (which holds almost surely under CJ$_{\beta, \delta}\sn$), then by
 \eqref{masterf}, $\Phi_{k,n}(1) \neq 0$, and 
one easily deduces that $G_k(U_n) \in E_k$, which allows to define $\log \Phi_{k,n}(1)$ without ambiguity. 
Now, the map from $\mathbb{D}^{n-1} \times (\mathbb{U} \backslash \{1\})$ to $\mathbb{R}$, given by 
$$(\gamma_1, \dots, \gamma_{n-1}) \mapsto \sum_{j=0}^{k-1} \log(1 -\gamma_j)$$
is continuous if we take the principal branch of the logarithm, 
and since $U_n$ depends continuously on $(\gamma_1, \dots, \gamma_{n-1}) \in \mathbb{D}^{n-1} \times (\mathbb{U} \backslash \{1\})$, as it can be checked in \cite{apa}, 
the map
$$(\gamma_1, \dots, \gamma_{n-1}) \mapsto \log \Phi_{k,n}(1)$$
is also continuous. 
These two maps have the same exponential, and one can check that they are both real if the $\gamma_j$'s are all real. Hence, they are equal, which fully justifies the equation
\begin{equation}
\label{sum}
\log \Phi_{k,n} (1) = \sum_{j=0}^{k-1} \log(1-\gamma_j).
\end{equation}

\subsection{Abel-Plana summation formula}
\begin{thm}
Let $m < n$ be integers and let $g$ be a holomorphic function on the strip
 $\{t \in \mathbb{C}, n \leq \Re\!\ t \leq m\}$ (i.e. $g$ is continuous on this strip and holomorphic in
its interior). We assume that
 $g(t) = o\left(\exp (2\pi |\Im\!\ t|)\right)$ as $\Im\!\ t \rightarrow \pm\infty$, uniformly with
 respect to $\Re\!\ t \in [n,m]$. Then, 
\begin{eqnarray}\nonumber
\sum_{j= m+1}^n g(j) &=&  \int_m^n g(t)dt + \frac{g(n) - g(m)}{2}\\
\nonumber
 &+&\ii \int_0^\infty \frac{g(m+\ii y) - g(n+ \ii y) - g(m -\ii y) + g(n-\ii y)}{e^{2\pi y}-1} dy\,. \\
\label{AP}
\end{eqnarray}
\end{thm}
For a proof see \cite{Olver} p. 290.

%\end{rem}
\bibliographystyle{plain}
\bibliography{a2j}

\end{document}